# INFERENCE FOR EIGENVALUES AND EIGENVECTORS OF GAUSSIAN SYMMETRIC MATRICES

By Armin Schwartzman,[1] Walter F. Mascarenhas[2]
and Jonathan E. Taylor[3]

*Harvard School of Public Health, Universidade de São Paulo
and Stanford University*

This article presents maximum likelihood estimators (MLEs) and log-likelihood ratio (LLR) tests for the eigenvalues and eigenvectors of Gaussian random symmetric matrices of arbitrary dimension, where the observations are independent repeated samples from one or two populations. These inference problems are relevant in the analysis of diffusion tensor imaging data and polarized cosmic background radiation data, where the observations are, respectively, $3 \times 3$ and $2 \times 2$ symmetric positive definite matrices. The parameter sets involved in the inference problems for eigenvalues and eigenvectors are subsets of Euclidean space that are either affine subspaces, embedded submanifolds that are invariant under orthogonal transformations or polyhedral convex cones. We show that for a class of sets that includes the ones considered in this paper, the MLEs of the mean parameter do not depend on the covariance parameters if and only if the covariance structure is orthogonally invariant. Closed-form expressions for the MLEs and the associated LLRs are derived for this covariance structure.

**1. Introduction.** Consider the signal-plus-noise model

$$(1) \qquad\qquad Y = M + Z,$$

where $Y, M, Z \in \mathcal{S}_p$, the set of $p \times p$ symmetric matrices ($p \geq 2$). Here $M$ (capital $\mu$) is a mean parameter and $Z$ is a mean-zero Gaussian random

Received June 2007; revised February 2008.

[1]Supported in part by a William R. and Sara Hart Kimball Stanford Graduate Fellowship and by the Statistical and Applied Mathematical Sciences Institute.

[2]Supported in part by CNPq projeto 309470/2006-4.

[3]Supported in part by NSF Grant DMS-04-05970.

matrix. Our goal is to estimate and test hypotheses about $M$ when it is restricted to subsets of $\mathcal{S}_p$ defined in terms of the eigenvalues and eigenvectors of $M$. In particular, we are interested in testing whether $M$ has a fixed set of eigenvalues or a fixed set of eigenvectors, taking into account ordering, and whether the eigenvalues of $M$ have particular multiplicities. We consider both the one- and two-sample problems, where repeated independent observations of $Y$ are sampled from one population with mean $M$ or two populations with means $M_1$ and $M_2$, respectively. In the two-sample problem, the pairs $(M_1, M_2)$ are restricted to subsets of $\mathcal{S}_p \times \mathcal{S}_p$ defined in terms of the eigenvalues and eigenvectors of $M_1$ and $M_2$. Here we are interested in testing whether $M_1$ and $M_2$ have the same eigenvalues or the same eigenvectors, taking into account the possibility of the eigenvalues having particular multiplicities.

This problem contains many of the ingredients of a classical multivariate problem. What makes it interesting and relevant are the following two aspects of it. First, new massive data sources have appeared recently where the observations take the form of random symmetric matrices, for which model (1) is appropriate. Second, the mean parameter sets involved in the inference problems for eigenvalues and eigenvectors have interesting geometries as subsets of $\mathcal{S}_p$ and $\mathcal{S}_p \times \mathcal{S}_p$: affine Euclidean subspaces, orthogonally invariant embedded submanifolds and convex polyhedral cones. Deriving maximum likelihood estimators (MLEs) in these sets is nontrivial as the multiplicities of the eigenvalues of $M$ affect both the dimension of the set and whether the true parameter falls on the boundary of the set. The effect of the covariance structure also depends on the geometry of the set. We show in Theorems 3.1 and 3.2 that for a particular class of sets that includes all the sets considered in this paper, the MLEs do not depend on the covariance parameters if and only if $Z$ in (1) has an orthogonally invariant covariance structure. Closed-form expressions for the MLEs and LLRs are then derived for this particular covariance structure in Propositions 4.1 and 4.2, and in Theorems 4.1, 4.2, 5.1 and their corollaries.

Examples of random symmetric matrix data are found in diffusion tensor imaging (DTI) and polarized cosmic background radiation (CMB) measurements. DTI is a modality of magnetic resonance imaging that produces at every voxel (3D pixel) a $3 \times 3$ symmetric positive definite matrix, also called a diffusion tensor (Basser and Pierpaoli [4], LeBihan et al. [20]). The diffusion tensor describes the local pattern of water diffusion in tissue. In the brain, it serves as a proxy for local anatomical structure. The tensor eigenvalues are generally indicative of the type of tissue and its health, while the eigenvectors indicate the spatial orientation of the underlying neural fibers. Since the anatomical information is contained in the eigenstructure of the diffusion tensor, inference methods for the eigenvalues and eigenvectors of the tensor are useful in anatomical imaging studies. The one- and two-sample testing problems arise when comparing images voxelwise to a fixed anatomical atlas or between two groups of subjects.

In astronomy, the polarization pattern of the CMB can be represented by $2 \times 2$ symmetric positive definite matrices; here again the information is contained in the eigenstructure. The eigenvalues relate to the electromagnetic field strength while the eigenvectors indicate the polarization orientation (Hu and White [15], Kogut et al. [17]). The one-sample estimation problem arises in image interpolation, while two-sample testing problems arise when comparing regions of the sky.

DTI data has been shown to be well modeled by (1) (Basser and Pajevic [3]). In this case the pattern of diffusion in the tissue is captured by $M$, while $Z$ captures the variability in the measurements. The Gaussianity assumption for $Z$ holds for high signal-to-noise ratio (SNR) and the orthogonal invariance property holds for orthogonally invariant field gradient designs. In the DTI literature, model (1) has also been used to model the data after a matrix log transformation (Arsigny et al. [1], Fletcher and Joshi [13], Schwartzman [30]). The matrix log, computed by taking the log of the eigenvalues and keeping the eigenvectors intact, maps the observed positive definite matrices to the set $\mathcal{S}_p$. This ensures that, when transformed back, the estimated matrices are always positive definite, as required by the anatomy. Whether the log transform should be applied is a subject of current debate (Whitcher et al. [33]). The methods developed in this paper are applicable in either case because the matrix log affects only the eigenvalues in a one-to-one fashion, so the various hypotheses about eigenvalues and eigenvectors can be equivalently stated in both domains.

The inference problems for the eigenvalues and eigenvectors of $M$ in (1) are nonstandard in that many of them involve nonlinear parameter sets. Suppose $M$ is constrained to lie in a generic subset $\mathcal{M} \subset \mathcal{S}_p$ or $\mathcal{M}_2 \subset \mathcal{S}_p \times \mathcal{S}_p$. The existence of MLEs and their consistency are guaranteed for closed parameter sets (Wald [32]), but their asymptotic distributions depend on the geometry of these sets. The cases considered here involve three kinds of sets:

- Affine subspaces of $\mathcal{S}_p$ and $\mathcal{S}_p \times \mathcal{S}_p$, that is, translations of a linear subspace by a constant $M_0 \in \mathcal{S}_p$ or $M_0 \in \mathcal{S}_p \times \mathcal{S}_p$. In these cases the MLE is unique and exactly multivariate normal.
- Embedded submanifolds of $\mathcal{S}_p$ and $\mathcal{S}_p \times \mathcal{S}_p$ defined by a set of constraints, much in the same way that the set $\mathcal{O}_p$ of $p \times p$ orthogonal matrices is defined by the set of constraints $Q'Q = QQ' = I_p$, where $I_p$ is the $p \times p$ identity matrix. If the submanifold is closed or if the true parameter is in the interior of the set, then the MLE is unique almost surely and asymptotically multivariate normal on the tangent space to the manifold at the true parameter (Efron [11]).

• Closed convex cones in $\mathcal{S}_p$, that is, closed convex sets $C$ where there exists a vertex $M_0$ such that $X \in C$ implies $a(X - M_0) + M_0 \in C$ for any real nonnegative $a$. In this case the MLE is also unique almost surely. The asymptotic distribution of the MLE is normal if the true parameter $M$ is inside the cone and is the projection of a normal onto the support cone at the true parameter if the true parameter is on the boundary of the cone (Self and Liang [31]).

In testing problems, the distribution of the LLRs associated with the MLEs depends on the geometries of the null set $\mathcal{M}_0$ and the alternative set $\mathcal{M}_A$. If the sets are nested and both are closed embedded submanifolds, then as $n \to \infty$, the LLR is asymptotically $\chi^2$ with number of degrees of freedom equal to the difference in dimension between the null and the alternative (e.g., Mardia, Kent and Bibby [23], page 124, Lehman [21], page 486). If the null set is a cone that is not a submanifold and the alternative is unrestricted, then the LLR is a mixture of $\chi^2$ with various numbers of degrees of freedom according to the dimensions of the faces of the cone and with weights that depend on the angles between the faces (Self and Liang [31]). The multiplicities of the eigenvalues of $M$ play a crucial role as they determine the dimension of the hypotheses. For example, if a LLR test statistic is derived assuming specific multiplicities but the true parameter has different multiplicities, then the LLR no longer has the prescribed asymptotic $\chi^2$ distribution (Drton [8]). For this reason, we assume throughout this paper that the multiplicities of the eigenvalues at any particular hypothesis are fixed and known. This is a reasonable assumption in DTI data, as different tissue types typically correspond to diffusion patterns with characteristic eigenvalue multiplicities, often called isotropic, prolate, oblate, and fully anisotropic diffusion (Zhu et al. [34]).

The normal distribution for symmetric random matrices has been known in statistics for over half a century. The orthogonally invariant covariance structure dates at least as far back as the work of Mallows [22], who presented it in its most general form. James [16] obtained a spherically symmetric version of this distribution as a particular limit of the Wishart distribution and applied it in a study of the ordered characteristic roots of random matrices. Another form of the normal distribution for symmetric matrices has been studied as a special case of the normal distribution for rectangular random matrices (Gupta and Nagar [14]). The covariance in this case is defined in terms of the rows and columns of the matrix and is not necessarily orthogonally invariant. A connection between this distribution and James' is made by Chikuse [7]. In modern random matrix theory, the most common Gaussian distribution for symmetric matrices is the Gaussian Orthogonal Ensemble (GOE), which was developed independently in the physics literature (Mehta [24]). The covariance structure of the GOE is orthogonally invariant and is a special case of the general orthogonally invariant covariance structure originally introduced by Mallows [22].

Despite the historical presence of the Gaussian distribution for symmetric matrices, little can be found in the statistics literature about inference for eigenvalues and eigenvectors when samples are drawn from that distribution, at least for fixed $p$. The only reference we have found is the work by Mallows [22], who studied linear hypotheses involved in testing the eigenvectors of a single matrix. Our expressions for the MLEs and the corresponding LLRs are derived assuming the same orthogonally invariant covariance structure described there. The results are stated in Propositions 4.1 and 4.2, and in Theorems 4.1, 4.2 and 5.1, and their corollaries. We show in Theorems 3.1 and 3.2 that the orthogonally invariant covariance structure is precisely the type of covariance needed to make the problem of finding MLEs independent of the covariance parameters for a class of parameter sets that include all the sets considered in the present article. Proposition 3.1 provides a LLR test for checking whether the orthogonally invariant covariance structure is appropriate for any particular data set.

The theory of curved exponential families is useful for finding the asymptotic distribution of the MLEs on embedded submanifolds. Interestingly, there is a close resemblance between the problem of estimating the mean of a normal symmetric matrix when its eigenvalues are fixed and Fisher's classical circle problem of a bivariate normal with mean parameter constrained to a circle of fixed radius (Efron [11]). Other useful tools are optimization techniques on $\mathcal{O}_p$ (Chang [5], Edelman, Arias and Smith [10]) and algorithms for estimation over parameter sets with edges (e.g., Chernoff [6], Self and Liang [31]), and sets with linear constraints (e.g., Lawson and Hanson [19], Dykstra [9]).

The inference problems considered in this paper are motivated by the analysis of DTI data. There are many combinations of inference problems for eigenvalues and eigenvectors depending on the restrictions that are imposed on the parameters. Section 2 provides an overview of the different cases. To consolidate the results, we first derive general results regarding the MLEs and LLRs for the orthogonally invariant symmetric-matrix-variate normal distribution in Section 3. We then apply these results to the various specific problems for eigenvalues and eigenvectors in Sections 4 and 5.

## 2. Inference for eigenvalues and eigenvectors: overview.

2.1. *One-sample problems.* The one-sample problem in DTI data is useful when comparing an individual or a small group of individuals to a fixed anatomical atlas. In this comparison, the investigator may be interested in making inferences about the full tensor $M$, only about its eigenvalues, or only about its eigenvectors. The cases considered here are summarized in

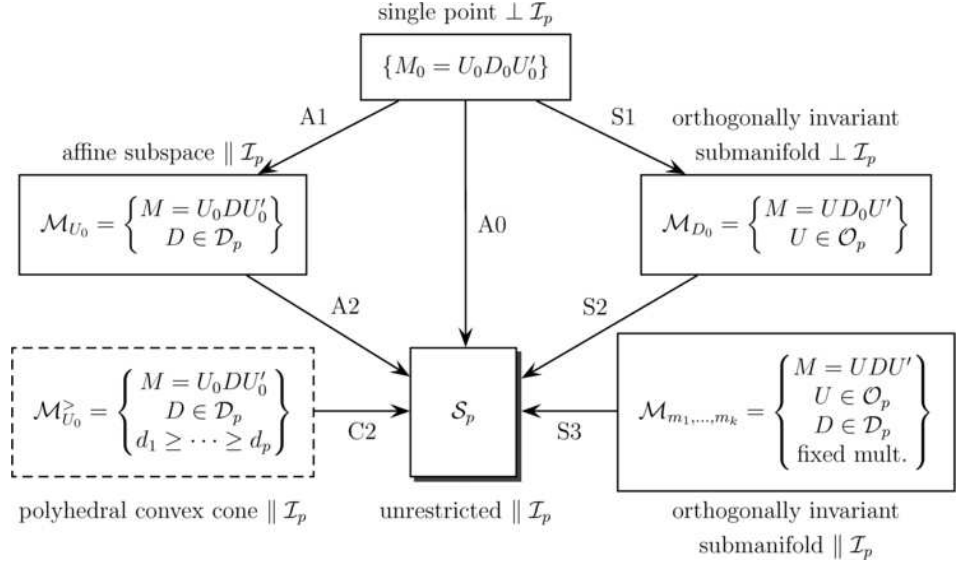

Fig. 1.   *Tests and parameter sets in the one-sample problem. Each test is represented by an arrow, where the origin of the arrow indicates the null hypothesis and the end of the arrow indicates the alternative hypothesis. The dashed borders indicate sets for which the MLE has no closed-form. The orthogonal and parallel notation is explained in Definition 3.1, Section 3.2.*

Figure 1. Of the seven tests, three involve affine subspaces, one involves a convex cone, and three involve orthogonally invariant submanifolds. All the sets involved satisfy the conditions of Theorem 3.1 below, so that the MLE of $M$ does not depend on the covariance parameters when the covariance structure is orthogonally invariant. We begin with the tests involving affine subspaces.

(A0) *Unrestricted test.* Here no particular attention is paid to the eigenstructure. The test is $H_0 : M = M_0$ vs. $H_A : M \neq M_0$.

(A1) *Test of eigenvalues with known unordered eigenvectors.* This test is useful for assessing differences in eigenvalues, while the eigenvectors are fixed from the atlas. Denote by $\mathcal{D}_p$ the set of $p \times p$ diagonal matrices. The test is $H_0 : M = M_0 = U_0 D_0 U_0'$ vs. $H_A : M \in \mathcal{M}_{U_0}$, where

$$\mathcal{M}_{U_0} = \{M = U_0 D U_0' : D \in \mathcal{D}_p\} \tag{2}$$

for fixed $U_0 \in \mathcal{O}_p$. This is equivalent to testing whether the eigenvalues $D = U_0' M U_0$ of $M$ are different in value or in order from the eigenvalues $D_0 = U_0' M_0 U_0$ of $M_0$. The fixed eigenvectors $U_0$ are assumed unordered in the sense that they represent only a coordinate system, thus allowing the eigenvalues to change their order under the alternative. Moreover, no restrictions are placed as to whether there are multiplicities in the eigenvalues of $D$. These assumptions make $\mathcal{M}_{U_0}$ a linear subspace of $\mathcal{S}_p$, a hyperplane of dimension $p$ that contains the multiples of $I_p$.

(A2) *Test of unordered eigenvectors with unknown eigenvalues.* This is another test that falls in the linear category. The test is $H_0 : M \in \mathcal{M}_{U_0}$ vs. $M \notin \mathcal{M}_{U_0}$, where $\mathcal{M}_{U_0}$ is given by (2). This is a test of whether $M$ has particular eigenvectors $U_0$ in any order (i.e., the columns of $U_0$ diagonalize $M$) while the eigenvalues are treated as nuisance parameters. This is one of the cases that was considered by Mallows [22].

Next we consider one test involving a polyhedral convex cone.

(C2) *Test of ordered eigenvectors with unknown eigenvalues.* This test serves the same purpose as test (A2), but here the eigenvectors are forced to be ordered, so that the first eigenvector corresponds to the largest eigenvalue, the second eigenvector to the second largest eigenvalue, and so on (ties allowed). The change in assumptions makes the problem nonlinear. Here the test is $H_0 : M \in \mathcal{M}_{U_0}^{>}$ vs. $M \notin \mathcal{M}_{U_0}^{>}$, where

$$(3) \qquad \mathcal{M}_{U_0}^{>} = \{M = U_0 D U_0' : D \in \mathcal{D}_p, d_1 \geq \cdots \geq d_p\},$$

$U_0 \in \mathcal{O}_p$ is fixed and $d_1 \geq \cdots \geq d_p$ are the diagonal entries of $D$. The parameter set $\mathcal{M}_{U_0}^{>}$ is a closed polyhedral convex cone determined by the constraints $d_1 \geq \cdots \geq d_p$ and rotated by the matrix $U_0$.

The last three tests involve orthogonally invariant embedded submanifolds of $\mathcal{S}_p$.

(S1) *Test of ordered eigenvectors with known eigenvalues.* This test serves the same purpose as test (A2), but here the eigenvalues are treated as fixed from the atlas. The change in assumptions makes the problem nonlinear. The test is $H_0 : M = M_0 = U_0 D_0 U_0'$ vs. $M \in \mathcal{M}_{D_0}$, where

$$(4) \qquad \mathcal{M}_{D_0} = \{M = U D_0 U' : U \in \mathcal{O}_p\}$$

for fixed $D_0 \in \mathcal{D}_p$ with diagonal entries $d_1 \geq \cdots \geq d_p$. If the diagonal entries of $D_0$ are distinct, this is equivalent to testing whether $M$ has eigenvector matrix $U_0 \in \mathcal{O}_p$ (up to sign flips of the columns of $U_0$), while the eigenvalues $d_1 > \cdots > d_p$ of $M$ are known and fixed. In general, the set $\mathcal{M}_{D_0}$ is the set of matrices with fixed $k \leq p$ distinct eigenvalues $\tilde{d}_1 > \cdots > \tilde{d}_k$ and corresponding multiplicities $m_1, \ldots, m_k$ so that $\sum_{i=1}^{k} m_i = p$. This is a closed embedded submanifold of $\mathcal{S}_p$ that is invariant under transformations $QMQ'$ for $Q \in \mathcal{O}_p$. Note that its dimension depends on the multiplicities $m_1, \ldots, m_k$. For example, if the eigenvalues have no multiplicities ($k = p$), then $\mathcal{M}_{D_0}$ is diffeomorphic to $\mathcal{O}_p$ and has dimension $p(p-1)/2$, but if all the eigenvalues are equal ($k = 1$, $m_1 =$ $p$), then $\mathcal{M}_{D_0}$ reduces to the single point $\{\tilde{d}_1 I_p\}$. In the case $p = 2$ with $d_1 \neq d_2$, the set $\mathcal{M}_{D_0}$ is exactly the same as the parameter set in Fisher's circle problem of a bivariate normal with mean parameter constrained to a circle of fixed radius (Efron [11]).

(S2) *Test of eigenvalues with unknown eigenvectors.* This test serves the same purpose as test (A1), except that the eigenvectors are treated as nuisance parameters. Again the change in assumptions makes the problem nonlinear. Here the test is $H_0 : M \in \mathcal{M}_{D_0}$ vs. $H_A : M \notin \mathcal{M}_{D_0}$, where $\mathcal{M}_{D_0}$ is given by (4). This is equivalent to testing whether $M$ has a set of eigenvalues equal to those in $D_0$, where the set may include repeats depending on the multiplicities. Because the eigenvectors are not specified, the null includes reorderings of the eigenvalues.

(S3) *Test of eigenvalue multiplicities.* This is a test of whether $M$ has eigenvalues with specific multiplicities against the eigenvalues being all distinct. This is useful in DTI to identify tissue types since, in theory, tensors in the cerebral fluid and the gray matter are isotropic ($d_1 = d_2 = d_3$), tensors in single fibers are prolate ($d_1 > d_2 = d_3$) and tensors in fiber crossings may be oblate ($d_1 = d_2 > d_3$). The test is $M \in \mathcal{M}_{m_1,\ldots,m_k}$ vs. $M \notin \mathcal{M}_{m_1,\ldots,m_k}$, where

$$(5) \qquad \mathcal{M}_{m_1,\ldots,m_k} = \{M = UDU' : U \in \mathcal{O}_p, D \in \mathcal{D}_p, \text{ mult. } m_1,\ldots,m_k\}$$

is the set of matrices with unspecified $k \leq p$ distinct eigenvalues $\tilde{d}_1 > \cdots > \tilde{d}_k$ and corresponding multiplicities $m_1,\ldots,m_k$ so that $\sum_{i=1}^k m_i = p$. This is an orthogonally invariant embedded submanifold of $\mathcal{S}_p$ whose dimension depends on the multiplicities. For example, if $d_1 = \cdots = d_p$ ($k = 1$, $m_1 = p$) we get the closed straight line $\mathcal{I}_p = \{\alpha I_p, \alpha \in \mathbb{R}\}$ of dimension 1, but if $d_1 \neq \cdots \neq d_p$ ($k = p$) then we get the open set $\mathcal{S}_p \setminus \mathcal{I}_p$ of dimension $p$.

The specific cases above are worked out in detail in Section 4. Before proceeding, it is helpful to visualize the above parameter sets in the case $p = 2$.

EXAMPLE 2.1 ($p = 2$). For $X \in \mathcal{S}_2$, define the embedding of $\mathcal{S}_2$ in $\mathbb{R}^3$ by the operator $(x, y, z) = \text{vecd}(X) = (X_{11}, X_{22}, \sqrt{2}X_{12})$. In Figure 2, the set of diagonal matrices $\text{vecd}(D) = (d_1, d_2, 0)$ corresponds to the $xy$-plane, while the set of matrices with equal eigenvalues $d_1 = d_2$ are multiples of the identity matrix and map to the line $\mathcal{I}_2 = \{x = y, z = 0\}$. For a generic $M$, consider the eigendecomposition

$$(6) \qquad M = UDU' = \begin{pmatrix} \cos\theta & -\sin\theta \\ \sin\theta & \cos\theta \end{pmatrix} \begin{pmatrix} d_1 & 0 \\ 0 & d_2 \end{pmatrix} \begin{pmatrix} \cos\theta & \sin\theta \\ -\sin\theta & \cos\theta \end{pmatrix}.$$

Then $M$ gets mapped to $(x, y, z) = \mathrm{vecd}(M)$, where

$$x = (d_1 + d_2)/2 + [(d_1 - d_2)/2]\cos 2\theta,$$
$$y = (d_1 + d_2)/2 - [(d_1 - d_2)/2]\cos 2\theta,$$
$$z = \sqrt{2}[(d_1 - d_2)/2]\sin 2\theta.$$

The effect of $U$ is a rotation around the axis $\mathcal{I}_2$ by an angle $2\theta$. In Figure 2, the circle represents the rotation trajectory of the point $(3, 1, 0)$ on the plane $x + y = 4$ with center at $(2, 2)$.

The set $\mathcal{M}_{U_0}$ defined by (2) has fixed $\theta$ and unrestricted $d_1$ and $d_2$, so it maps to a plane through the line $\mathcal{I}_2$ making an angle $2\theta$ with the $xy$ plane. Notice that replacing $\theta$ by $\theta + \pi$ results in the same plane, as the eigendecomposition is invariant under sign changes of the eigenvectors. If the eigenvalues are constrained so that $d_1 \geq d_2$, then we get the set $\mathcal{M}_{U_0}^{>}$ defined by (3). This is half the previous plane, a closed convex cone that is not an embedded submanifold of $\mathcal{S}_2$. The set $\mathcal{M}_{D_0}$ defined by (4) has fixed $d_1$ and $d_2$ and allows $\theta$ to vary. When $d_1 \neq d_2$ it maps to a circle of radius $\sqrt{2}(d_1 - d_2)/2$ orthogonal to the line $\mathcal{I}_2$ and centered at $((d_1 + d_2)/2, (d_1 + d_2)/2, 0)$, shown in Figure 2 for $d_1 = 3$ and $d_2 = 1$. This is a closed embedded submanifold of $\mathcal{S}_2$ that is invariant under orthogonal transformations $QMQ'$ for $Q \in O(2)$. When projected on the plane of the circle, this parameter set is exactly the same as the parameter set in Fisher's circle problem of a bivariate normal with mean parameter constrained to a circle of fixed radius (Efron [11]). The dimension of this set depends on the multiplicities of the eigenvalues. If $d_1 = d_2$, then the circle collapses to a single point.

2.2. *Two-sample problems.* The two-sample problem in DTI data is useful when comparing two groups, such as in case-control studies. As in the one-sample case, there are many combinations of hypotheses. Here we consider nontrivial situations where some parameters are common and others are not: the case where $M_1$ and $M_2$ have common eigenvalues and the case where $M_1$ and $M_2$ have common eigenvectors. The cases considered here are summarized in Figure 3. Of the five tests, only the first is linear. The rest are nonlinear. All the sets involved satisfy the conditions of Theorem 3.2 below, so that the MLEs of $M_1$ and $M_2$ do not depend on the covariance parameters when the covariance structure is orthogonally invariant.

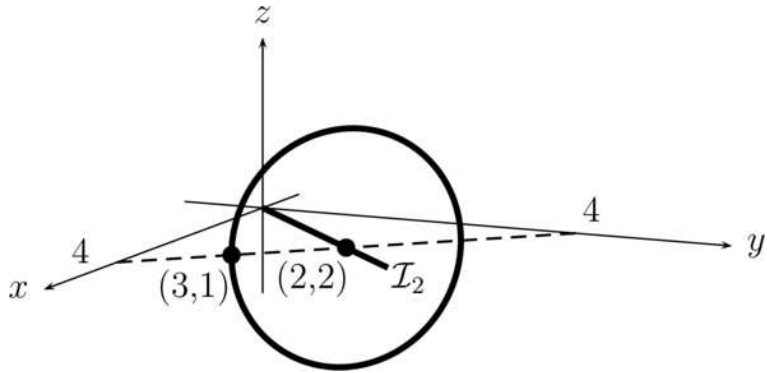

Fig. 2. *Examples of orthogonally invariant subsets of $\mathcal{S}_2$ embedded in $\mathbb{R}^3$.*

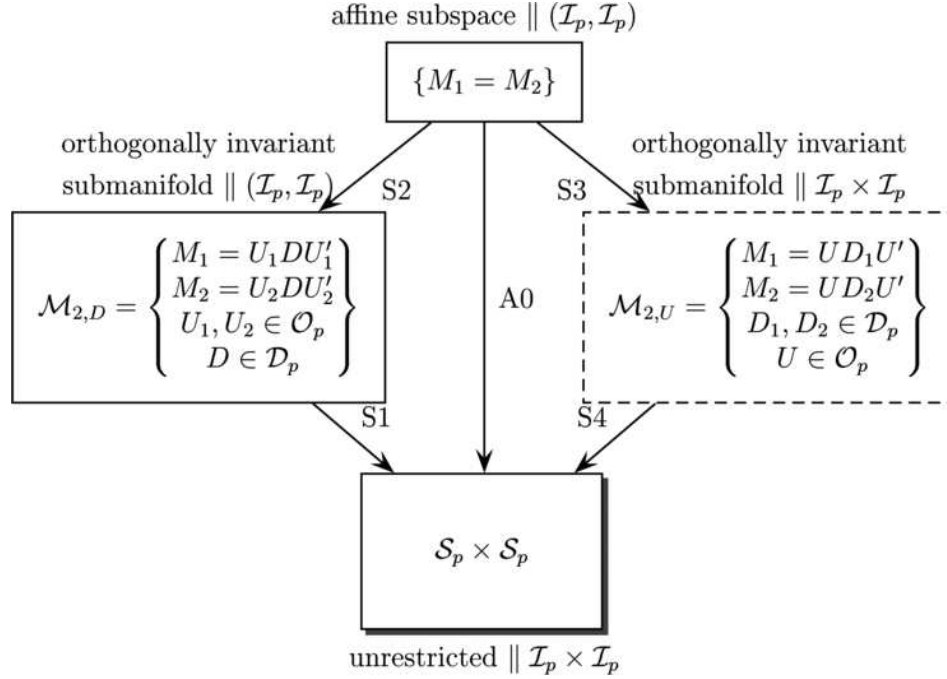

Fig. 3.   *Test and parameter sets in the two-sample problem. Each test is represented by an arrow, where the origin of the arrow indicates the null hypothesis and the end of the arrow indicates the alternative hypothesis. The dashed borders indicate sets for which the MLE has no closed-form. The orthogonal and parallel notation is explained in Definition 3.2, Section 3.3.*

(A0)  *Unrestricted test.* Here no particular attention is paid to the eigenstructure. The test is $H_0 : M_1 = M_2$ vs. $H_A : M_1 \neq M_2$.

(S1)  *Test of equality of eigenvalues with unrestricted eigenvectors.* This is a test of whether the means of two populations have the same eigenvalues, treating the eigenvectors as nuisance parameters. In DTI, this is useful for determining, for instance, if the two populations correspond to the same tissue type. The test is $H_0 : (M_1, M_2) \in \mathcal{M}_{2,D}$ vs.

$H_A \colon (M_1, M_2) \notin \mathcal{M}_{2,D}$, where

$$(7) \qquad \mathcal{M}_{2,D} = \{(M_1, M_2) \colon M_1 = U_1 D U_1', M_2 = U_2 D U_2'\}$$

for unspecified $D \in \mathcal{D}_p$ and $U_1, U_2 \in \mathcal{O}_p$. $D$ is assumed to have unspecified $k$ distinct eigenvalues $\tilde{d}_1 > \cdots > \tilde{d}_k$ with fixed and known multiplicities $m_1, \ldots, m_k$. This is an embedded submanifold that is invariant under transformations $Q_1 M_1 Q_1'$ and $Q_2 M_2 Q_2'$ for $Q_1, Q_2 \in \mathcal{O}_p$. Its dimension depends on the multiplicities.

(S2) *Test of equality of eigenvectors with common eigenvalues.* This is a test of whether the means of two populations have the same eigenvectors, when the eigenvalues are treated as nuisance parameters and assumed equal between the two populations. In DTI, this is useful for determining, for instance, if two populations of the same tissue type have neural fibers oriented in the same direction in space. The test is $H_0 \colon M_1 = M_2$ vs. $H_A \colon (M_1, M_2) \in \mathcal{M}_{2,D}$, where $\mathcal{M}_{2,D}$ is given by (7). Here again the dimension depends on the multiplicities of the eigenvalues in $D$, which are assumed to be known.

(S3) *Test of equality of eigenvalues with common eigenvectors.* This test serves the same purpose as test (S1) of testing equality of eigenvalues while treating the eigenvectors as nuisance parameters, with the added restriction that the eigenvectors are common to both populations. The test is $H_0 \colon M_1 = M_2$ vs. $H_A \colon (M_1, M_2) \in \mathcal{M}_{2,U}$, where

$$(8) \qquad \mathcal{M}_{2,U} = \{(M_1, M_2) \colon M_1 = U D_1 U', M_2 = U D_2 U'\}$$

for unspecified $D_1, D_2 \in \mathcal{D}_p$ and $U \in \mathcal{O}_p$. Again, this is an embedded submanifold that is invariant under transformations $Q M_1 Q'$ and $Q M_2 Q'$ for $Q \in \mathcal{O}_p$. Its dimension depends on the multiplicities of the eigenvalues in $D_1$ and $D_2$.

(S4) *Test of equality of eigenvectors with unrestricted eigenvalues.* This test serves the same purpose as test (S2) of testing equality of eigenvectors, when the eigenvalues are treated as nuisance parameters, but the assumption of common eigenvalues has been removed. Here the test is $H_0 \colon (M_1, M_2) \in \mathcal{M}_{2,U}$ vs. $H_A \colon (M_1, M_2) \notin \mathcal{M}_{2,U}$, where $\mathcal{M}_{2,U}$ is given by (8).

The set $\mathcal{M}_{2,U}$ involved in cases (S3) and (S4) is a difficult case for which the MLE does not have a closed-form. For this reason, below we focus on the solutions to cases (S1) and (S2). These are worked out in detail in Section 5.

**3. General inference for Gaussian symmetric matrices.**

3.1. *Covariance structures.* Consider the embedding of $\mathcal{S}_p$ in $\mathbb{R}^q$, $q = p(p+1)/2$, $p \geq 2$, by the operator $\text{vecd}(Y) = (\text{diag}(Y)', \sqrt{2} \, \text{offdiag}(Y)')'$, where $\text{diag}(Y)$ is a $p \times 1$ column vector containing the diagonal entries of $Y \in \mathcal{S}_p$ and $\text{offdiag}(Y)$ is a $(q-p) \times 1$ column vector containing the off-diagonal entries of $Y$. A random matrix $Y \in \mathcal{S}_p$ is nondegenerate Gaussian if and only if it has density

$$(9) \qquad f(Y) = \frac{1}{(2\pi)^{q/2}|\Sigma|^{1/2}} \exp\left(-\frac{1}{2}\|Y - M\|_\Sigma^2\right)$$

with respect to Lebesgue measure on $\mathbb{R}^q$, where $\Sigma \in \mathcal{S}_q$ is positive definite and $\|\cdot\|_\Sigma^2$ is the norm corresponding to the inner product in $\mathcal{S}_p$

$$(10) \qquad \langle A, B \rangle_\Sigma = (\text{vecd}(A))'\Sigma^{-1}\text{vecd}(B).$$

The density (9), which we denote $N_{pp}(M, \Sigma)$, has mean parameter $M = E(Y) \in \mathcal{S}_p$ and covariance parameter $\Sigma = \text{cov}(\text{vecd}(Y))$. Notice that $\|Y - M\|_\Sigma^2 \sim \chi^2(q)$.

The covariance structure of $Z = Y - M$ in (9) is called orthogonally invariant if the distribution of $Z$ is the same as that of $QZQ'$ for any $Q \in \mathcal{O}_p$. Let $\mathbf{1}_p = (1, \ldots, 1)'$ be the $p$-vector of ones. Mallows [22] showed that $Z$ has orthogonally invariant covariance if and only if $\text{cov}(\text{diag}(Z)) = \sigma^2(I_p + c\mathbf{1}_p\mathbf{1}_p')$ for some $\sigma^2$ and $c$, and $\text{cov}(\text{offdiag}(Z)) = (\sigma^2/2)I_{q-p}$, independent of $\text{diag}(Z)$. For $p = 3$, the orthogonally invariant covariance $\Sigma$ is

$$\Sigma = \sigma^2 \begin{pmatrix} 1+c & c & c & 0 & 0 & 0 \\ c & 1+c & c & 0 & 0 & 0 \\ c & c & 1+c & 0 & 0 & 0 \\ 0 & 0 & 0 & 1 & 0 & 0 \\ 0 & 0 & 0 & 0 & 1 & 0 \\ 0 & 0 & 0 & 0 & 0 & 1 \end{pmatrix}.$$

Positive definiteness of $\Sigma$ requires that $\sigma^2 > 0$ and $c > -1/p$. In particular when $c \geq 0$, $Z$ may be constructed as $Z = \sigma(\sqrt{c}I_p w + W)$, where $w \sim N(0, 1)$ and $W \in \mathcal{S}_p$ has the GOE distribution, that is, $W$ has independent diagonal entries $N(0, 1)$ and independent off-diagonal entries $N(0, 1/2)$.

For convenience, we define $\tau = c/(1 + pc) < 1/p$. The distribution of $Y = M + Z$, which we call the orthogonally invariant symmetric-matrix-variate normal distribution, and denote $N_{pp}(M, \sigma^2, \tau)$, has density

$$(11) \qquad f(Y) = \frac{\sqrt{1 - p\tau}}{(2\pi)^{q/2}\sigma^q} \exp\left(-\frac{1}{2}\|Y - M\|_{\sigma^2, \tau}^2\right),$$

where the inner product (10) defining the norm $\|\cdot\|_{\sigma^2, \tau}^2$ has the explicit form

$$(12) \qquad \langle A, B \rangle_{\sigma^2, \tau} = [\text{tr}(AB) - \tau \, \text{tr}(A) \, \text{tr}(B)]/\sigma^2.$$

Notice that $\|A\|_{1,0}^2 = \operatorname{tr}(A^2)$ is the standard Frobenius norm of symmetric matrices. In the density (11), the mean parameter is $M \in \mathcal{S}_p$, while the covariance is captured by the two scalar parameters $\sigma^2 > 0$ and $\tau \in (-\infty, 1/p)$. In this notation, the GOE distribution is the same as $N_{pp}(0,1,0)$ and may be thought of as a standard normal for symmetric matrices. In general when $\tau = 0$, we call the covariance structure spherical. Later we use the following two obvious properties derived from the orthogonal invariance:

(1) If $Y \sim N_{pp}(M, \sigma^2, \tau)$, then $QYQ' \sim N_{pp}(QMQ', \sigma^2, \tau)$ for all $Q \in \mathcal{O}_p$.
(2) $\langle A, B \rangle_{\sigma^2, \tau} = \langle QAQ', QBQ' \rangle_{\sigma^2, \tau}$ for all $Q \in \mathcal{O}_p$.

Recall from Example 2.1 that when $p = 2$, the effect of the operation $QZQ'$ for $Q \in \mathcal{O}(2)$ is a rotation of $\operatorname{vecd}(Z)$ around the axis $\mathcal{I}_2$. The orthogonal invariance of the distribution of $Z$ implies that the contours of the density $N_{22}(M, \sigma^2, \tau)$ are ellipsoids that are circularly symmetric around an axis that is parallel to $\mathcal{I}_2$. Similarly for general $p$, the contours of (11) are hyper-ellipsoids that are spherically symmetric around an axis that is parallel to the line $\mathcal{I}_p = \{\alpha I_p, \alpha \in \mathbb{R}\}$. The cross sections to that axis are hyperspheres of dimension $q - 1$. The parameter $\tau$ only affects the size of the unequal axis, where the ratio between the length of the unequal axis and the diameter of the hyperspheres is $\sqrt{1 + pc} = 1/\sqrt{1 - p\tau}$. An important consequence is that $\tau$ only affects the variability parallel to the $\mathcal{I}_p$ axis, whereas the variability in the space orthogonal to $\mathcal{I}_p$ is the same as if the covariance structure were spherical. The orthogonally invariant covariance is the only one with this property. This is stated formally in Lemma 3.1 below. Notice that, by (12), a vector $A \in \mathcal{S}_p$ is orthogonal to $I_p$, that is, $\langle A, I_p \rangle_{\sigma^2, \tau} = 0$, if and only if $\operatorname{tr}(A) = 0$.

LEMMA 3.1. $\Sigma$ is orthogonally invariant if and only if there exists $\sigma^2$ such that, for all matrices $A, B \in \mathcal{S}_p$ with $\operatorname{tr}(A) = 0$, $\langle A, B \rangle_\Sigma = \langle A, B \rangle_{\sigma^2, 0}$.

PROOF. (i) Suppose $\Sigma$ is orthogonally invariant, so that $\langle \cdot, \cdot \rangle_\Sigma = \langle \cdot, \cdot \rangle_{\tilde{\sigma}^2, \tau}$ for some $\tilde{\sigma}^2, \tau$. Let $A$ be any matrix such that $\operatorname{tr}(A) = 0$. For any $B \in \mathcal{S}_p$,

$$\langle A, B \rangle_{\tilde{\sigma}^2, \tau} = [\operatorname{tr}(AB) - \tau \operatorname{tr}(A) \operatorname{tr}(B)]/\tilde{\sigma}^2 = \operatorname{tr}(AB)/\tilde{\sigma}^2 = \langle A, B \rangle_{\tilde{\sigma}^2, 0}$$

and we can take $\sigma^2 = \tilde{\sigma}^2$.

(ii) Conversely, suppose that there exists $\sigma^2$ such that $\langle A, B \rangle_\Sigma = \langle A, B \rangle_{\sigma^2, 0}$ for any pair of matrices $A, B$ such that $\operatorname{tr}(A) = 0$. Define the bilinear form $h(X, W) = \langle X, W \rangle_\Sigma - \langle X, W \rangle_{\sigma^2, 0}$. By (12), it is sufficient to show that $h(X, W) = \tau \operatorname{tr}(X) \operatorname{tr}(W)/\sigma^2$ for some $\tau$. By assumption, $h(X, W) = 0$ whenever $\operatorname{tr}(X) = 0$ or $\operatorname{tr}(W) = 0$. Thus, by bilinearity,

$$h(X, W) = h[X - \operatorname{tr}(X)I_p/p + \operatorname{tr}(X)I_p/p, W - \operatorname{tr}(W)I_p/p + \operatorname{tr}(W)I_p/p]$$

$$= h[X - \operatorname{tr}(X)I_p/p, W - \operatorname{tr}(W)I_p/p] + \operatorname{tr}(X)h[I_p/p, W - \operatorname{tr}(W)I_p/p]$$
$$\quad + \operatorname{tr}(W)h[X - \operatorname{tr}(X)I_p/p, I_p/p] + \operatorname{tr}(X)\operatorname{tr}(W)h(I_p/p, I_p/p)$$
$$= \operatorname{tr}(X)\operatorname{tr}(W)h(I_p, I_p)/p^2.$$

Therefore, $\Sigma$ is orthogonally invariant with parameters $\sigma^2$ and $\tau = \sigma^2 h(I_p, I_p)/p^2 = \sigma^2 \|I_p\|_\Sigma^2/p^2 - 1/p$. $\quad \square$

3.2. *General inference in the one-sample problem.* Let $Y_1, \ldots, Y_n$ i.i.d. $N_{pp}(M, \Sigma)$. Given $\Sigma$, $\bar{Y} = \sum_{i=1}^n Y_i/n$ is a sufficient and complete statistic for $M$. For future reference, the classical recipe for finding the MLE of $M$ over generic subsets of $\mathcal{S}_p$ is stated in the following lemma.

LEMMA 3.2. *Suppose $M$ lies in a nonempty closed subset $\mathcal{M}$ of $\mathcal{S}_p$. Given $\Sigma$, the MLE of $M$ over $\mathcal{M}$, denoted $\hat{M}$, minimizes the squared Mahalanobis distance*

$$(13) \qquad g_\Sigma(M) = \|\bar{Y} - M\|_\Sigma^2$$

*over $\mathcal{M}$ and is unique almost surely. $\hat{M}$ is consistent as $n \to \infty$.*

PROOF. Ignoring the constant $2\pi$, the log-likelihood from model (9) can be written as

$$(14) \qquad l(M, \Sigma) = -\frac{n}{2}\log|\Sigma| - \frac{n}{2}(s_\Sigma^2 + \|\bar{Y} - M\|_\Sigma^2),$$

where $s_\Sigma^2 = (1/n)\sum_{i=1}^n \|Y_i - \bar{Y}\|_\Sigma^2$ and we have used the sum-of-squares decomposition

$$\sum_{i=1}^n \|Y_i - M\|_\Sigma^2 = \sum_{i=1}^n \|Y_i - \bar{Y}\|_\Sigma^2 + n\|\bar{Y} - M\|_\Sigma^2.$$

Thus, given $\Sigma$, maximizing the likelihood over $\mathcal{M}$ is equivalent to minimizing (13) over $\mathcal{M}$. The existence, almost sure uniqueness and consistency of the MLE is guaranteed under the assumption that the set $\mathcal{M}$ is closed (Wald [32]). $\quad \square$

Lemma 3.2 says that $\hat{M}$ can be found as the orthogonal projection of $\bar{Y}$ on $\mathcal{M}$ according to the inner product $\langle A, B \rangle_\Sigma$. In general, $\hat{M}$ depends on $\Sigma$. However, if $\Sigma$ is orthogonally invariant, then for a particular class of sets $\mathcal{M}$ given by Definition 3.1 below, the projection $\hat{M}$ does not depend on $\Sigma$ and the problem is equivalent to finding the MLE when the covariance structure is spherical. This property is instrumental in the derivation of closed-form MLEs for the eigenvalue and eigenvector problems of Sections 4 and 5. Conversely, if we require the irrelevance of $\Sigma$ for finding the MLE of $M$ in that class of sets, then $\Sigma$ must be orthogonally invariant. This result is stated in Theorem 3.1 below.

DEFINITION 3.1. Define the straight line $\mathcal{I}_p = \{\alpha I_p, \alpha \in \mathbb{R}\}$ and let "$\oplus$" denote Minkowski addition. We say that a set $\mathcal{M}$ is:

(a) Orthogonal to $\mathcal{I}_p$, if it can be written as $\mathcal{M} = \mathcal{A} \oplus \{a I_p\}$ for some $a \in \mathbb{R}$ and some $\mathcal{A} \subset \mathcal{S}_p$ such that $\mathrm{tr}(A) = 0$ for all $A \in \mathcal{A}$.

(b) Parallel to $\mathcal{I}_p$, if it can be written as $\mathcal{M} = \mathcal{A} \oplus \mathcal{I}_p$, where $\mathcal{A}$ is orthogonal to $\mathcal{I}_p$.

The definition of orthogonality above is equivalent to the condition that for every tangent vector $A$ to $\mathcal{M}$, $\langle A, I_p \rangle_{\sigma^2, \tau} = 0$. For an orthogonal set $\mathcal{M}$, $\mathrm{tr}(M) = ap$ is constant for all $M \in \mathcal{M}$. Both properties defined in Definition 3.1 are orthogonally invariant: if $\mathcal{M}$ is orthogonal (parallel) to $\mathcal{I}_p$, so is the set $Q\mathcal{M}Q' = \{QMQ', M \in \mathcal{M}\}$, where $Q \in \mathcal{O}_p$. It is easy to check that each of the sets in Figure 1 belongs to one of the above categories: the single point $\{M_0\}$ and the set $\mathcal{M}_{D_0}$ are orthogonal to $\mathcal{I}_p$; the rest are parallel to $\mathcal{I}_p$.

THEOREM 3.1. $\Sigma$ is orthogonally invariant if and only if, for every embedded submanifold $\mathcal{M} \subset \mathcal{S}_p$ that is either orthogonal or parallel to $\mathcal{I}_p$, every critical point $\hat{M}$ of $g_{\Sigma}(M) = \|\bar{Y} - M\|_{\Sigma}^2$ over $\mathcal{M}$ is also a critical point of $g_{1,0}(M) = \|\bar{Y} - M\|_{1,0}^2 = \mathrm{tr}[(\bar{Y} - M)^2]$.

PROOF. (i) Suppose $\Sigma$ is orthogonally invariant, so $\langle \cdot, \cdot \rangle_{\Sigma} = \langle \cdot, \cdot \rangle_{\sigma^2, \tau}$ for some $\sigma^2, \tau$. Let $\mathcal{A}$ denote the tangent space to $\mathcal{M}$ at $\hat{M}$. For every critical point $\hat{M}$ of $g_{\sigma^2, \tau}(M) = \|\bar{Y} - M\|_{\sigma^2, \tau}^2$, $B = \bar{Y} - \hat{M}$ is orthogonal to all tangent vectors $A \in \mathcal{A}$, that is, $\langle A, B \rangle_{\sigma^2, \tau} = 0$. We need to show that every tangent vector $A$ also satisfies the spherical orthogonality condition $\langle A, B \rangle_{\sigma^2, 0} = \mathrm{tr}(AB) = 0$ for every $\tau$. By (12), this is achieved if either $\mathrm{tr}(A) = 0$ or $\mathrm{tr}(B) = 0$.

If $\mathcal{M}$ is orthogonal to $\mathcal{I}_p$, then $\mathrm{tr}(A) = 0$, and the result follows from Lemma 3.1. Assume instead $\mathcal{M}$ is parallel to $\mathcal{I}_p$. Decompose every tangent vector $A \in \mathcal{A}$ as a sum of a multiple of the identity $\mathrm{tr}(A) I_p / p$ and an orthogonal tangent vector $A - \mathrm{tr}(A) I_p / p$. Lemma 3.1 takes care of the orthogonal vectors. As for the remaining direction parallel to $\mathcal{I}_p$, note that $0 = \langle I_p, B \rangle_{\sigma^2, \tau} = \mathrm{tr}(B)(1 - p\tau)/\sigma^2$, which implies that $\mathrm{tr}(B) = 0$ for any critical point $\hat{M}$ for the inner product $\langle \cdot, \cdot \rangle_{\sigma^2, \tau}$. Therefore, $\langle I_p, B \rangle_{\sigma^2, 0} = \mathrm{tr}(B)/\sigma^2 = 0$.

(ii) Conversely, suppose now that for every $\bar{Y} \in \mathcal{S}_p$ and critical point $\hat{M} \in \mathcal{M}$ with respect to the inner product $\langle \cdot, \cdot \rangle_{\Sigma}$, $\hat{M}$ is also a critical point with respect to the inner product $\langle \cdot, \cdot \rangle_{1,0}$. Further, this holds for all sets that are orthogonal or parallel to $\mathcal{I}_p$. In particular, it holds for every orthogonal affine subspace of $\mathcal{S}_p$. Following a simple regression argument this implies that if we restrict our attention to the largest orthogonal subspace

$$\mathcal{A}_{\max} = \{A \in \mathcal{S}_p : \mathrm{tr}(A) = 0\},$$

then $\langle \cdot, \cdot \rangle_\Sigma$ agrees with $\langle \cdot, \cdot \rangle_{\sigma^2,0}$ for some $\sigma^2$.

Consider the case when the minimization problem is projection onto the subspace $\mathcal{A}_{\max}$. As the problem is simply a linear regression, there is a unique critical point for either $\langle \cdot, \cdot \rangle_\Sigma$ or $\langle \cdot, \cdot \rangle_{1,0}$ and, further, the critical points are identical by assumption. It is not hard to see that the critical point is $\hat{M} = \bar{Y} - \bar{Y} I_p/p$. The normal equations for this problem are

$$0 = \langle \bar{Y} - \hat{M}, A \rangle_\Sigma = \frac{\operatorname{tr}(\bar{Y})}{p} \langle I, A \rangle_\Sigma = \frac{\operatorname{tr}(\bar{Y})}{p} \langle I, A \rangle_{1,0}$$

and these must hold for all $\bar{Y}$ and all $A \in \mathcal{A}_{\max}$. In other words, $\langle I_p, A \rangle_{1,0} = \langle I_p, A \rangle_\Sigma = 0$, for all $A \in \mathcal{A}_{\max}$. Define $h(X, W) = \langle X, W \rangle_\Sigma - \langle X, W \rangle_{\sigma^2,0}$ where $\sigma^2$ is described above. We have to show that $h(X, Y) = \tau \operatorname{tr}(X) \operatorname{tr}(W)$ for some $\tau$. The proof proceeds similarly to the final display in the proof of Lemma 3.1. $\square$

By Theorem 3.1, if we assume the orthogonally invariant model (11), then within the class of sets $\mathcal{M}$ that are either orthogonal or parallel to $\mathcal{I}_p$, the problem of finding the MLE of $M$ reduces to minimizing $\operatorname{tr}[(\bar{Y} - M)^2]$ over $\mathcal{M}$. Once $\hat{M}$ is found, the estimates of $\sigma$ and $\tau$ in model (11) easily follow and are given by the following lemma.

LEMMA 3.3.  *Let $Y_1, \ldots, Y_n$ i.i.d. $N_{pp}(M, \sigma^2, \tau)$, $p \geq 2$. Suppose $M$ lies in a nonempty closed subset $\mathcal{M}$ of $\mathcal{S}_p$ and denote by $\hat{M}$ the MLE of $M$. Given $\tau$, if we define*

$$(15) \qquad \hat{\sigma}_\tau^2 = s_\tau^2 + \frac{1}{q} \|\bar{Y} - \hat{M}\|_{1,\tau}^2, \qquad s_\tau^2 = \frac{1}{qn} \sum_{i=1}^n \|Y_i - \bar{Y}\|_{1,\tau}^2,$$

*then $(\hat{M}, \hat{\sigma}_\tau^2)$ is the MLE of $(M, \sigma^2)$ in $\mathcal{M} \times \mathbb{R}_+$. Moreover, if we define*

$$(16) \qquad \hat{\tau} = -\frac{\sum_{i=1}^n \|Y_i - \bar{Y}\|_{1,q/p}^2 + n\|\bar{Y} - \hat{M}\|_{1,q/p}^2}{(q-1)[\sum_{i=1}^n [\operatorname{tr}(Y_i - \bar{Y})]^2 + n[\operatorname{tr}(\bar{Y} - \hat{M})]^2]},$$

*then $(\hat{M}, \hat{\sigma}_\tau^2, \hat{\tau})$ is the MLE of $(M, \sigma^2, \tau)$ in $\mathcal{M} \times \mathbb{R}_+ \times (-\infty, 1/p)$. The above estimates are unique almost surely and are consistent as $n \to \infty$.*

PROOF.  Under model (11), the log-likelihood (14) becomes

$$(17) \quad l(M, \sigma^2, \tau) = -\frac{nq}{2} \log \sigma^2 + \frac{n}{2} \log(1 - p\tau) - \frac{nq}{2\sigma^2} \left( s_\tau^2 + \frac{1}{q} \|\bar{Y} - M\|_{1,\tau}^2 \right).$$

The MLEs for $\sigma^2$ and $\tau$ are obtained by differentiating (17) and setting the derivative equal to zero. It is easy to check that $\hat{\tau} < 1/p$. Since $\hat{M} \to M$ as

$n \to \infty$,

$$\hat{\sigma}_\tau^2 = \frac{1}{qn}\sum_{i=1}^{n}\|Y_i - \hat{M}\|_{1,\tau}^2 \longrightarrow \frac{1}{q}E\|Y - M\|_{1,\tau}^2 = \sigma^2.$$

Similarly,

$$\hat{\tau} = -\frac{(1/n)\sum_{i=1}^{n}\|Y_i - \hat{M}\|_{1,q/p}^2}{((q-1)/n)\sum_{i=1}^{n}[\operatorname{tr}(Y_i - \hat{M})]^2} \longrightarrow -\frac{E\|Y - M\|_{1,q/p}^2}{(q-1)E[\operatorname{tr}(Y - M)]^2}.$$

The numerator is $E\|Y - M\|_{1,q/p}^2 = \operatorname{tr}[(I_p + \mathbf{1}\mathbf{1}'c)(I_p - \mathbf{1}\mathbf{1}'q/p)]\sigma^2 + (q-p)\sigma^2 = (1-q)pc\sigma^2$, where $c = \tau/(1-p\tau)$. The denominator is $(q-1)E[\operatorname{tr}(Y-M)]^2 = (q-1)\mathbf{1}'(I_p + \mathbf{1}\mathbf{1}'c)\mathbf{1}\sigma^2 = (q-1)(p + p^2c)\sigma^2$. Therefore $\hat{\tau} \to c/(1+pc) = \tau$. $\quad\square$

The distribution of the variance estimate (15) for known $\tau$ is established by noting that $n\|\bar{Y} - \hat{M}\|_{\sigma^2,\tau}^2 \to \chi^2(q-k)$ as $n \to \infty$ when $\mathcal{M}$ is a closed embedded submanifold of $\mathcal{S}_p$ of dimension $k$, and exactly so for all finite $n$ if $\mathcal{M}$ is an affine subspace. Multiplying (15) by the appropriate constants gives that, given $\tau$, $qn\hat{\sigma}_\tau^2/\sigma^2$ is the sum of two independent terms and thus has an exact or approximate $\chi^2$ distribution with $q(n-1) + q - k = qn - k$ degrees of freedom, depending on whether the distribution of the second term is exact or approximate for large $n$.

In the testing situation, let $Y_1, \ldots, Y_n \in \mathcal{S}_p$ i.i.d. $N_{pp}(M, \Sigma)$ and consider the test $H_0 : M \in \mathcal{M}_0$ vs. $H_A : M \in \mathcal{M}_A$, where $\mathcal{M}_0 \subset \mathcal{M}_A \subseteq \mathcal{S}_p$ are closed embedded submanifolds with dimensions $k_A = \dim(\mathcal{M}_A)$, $k_0 = \dim(\mathcal{M}_0)$, $k_A > k_0$. Let $\hat{M}_0$ and $\hat{M}_A$ denote MLEs of $M$ under $H_0$ and $H_A$, respectively, with corresponding maximized log-likelihoods $l_0$ and $l_A$.

LEMMA 3.4. *Suppose $\mathcal{M}_0 \subset \mathcal{M}_A$ are closed embedded submanifolds of $\mathcal{S}_p$. Given $\Sigma$, the LLR test statistic $2(l_A - l_0)$ is*

(18a) $$T = n\|\bar{Y} - \hat{M}_0\|_\Sigma^2 - n\|\bar{Y} - \hat{M}_A\|_\Sigma^2$$

(18b) $$= 2n\langle\bar{Y}, \hat{M}_A - \hat{M}_0\rangle_\Sigma + n\|\hat{M}_0\|_\Sigma^2 - n\|\hat{M}_A\|_\Sigma^2$$

*and is distributed as $\chi^2(k_A - k_0)$ under $M_0$ asymptotically as $n \to \infty$. If $\mathcal{M}_0 \subset \mathcal{M}_A$ are both affine subspaces of $\mathcal{S}_p$, then the distribution is exactly $\chi^2(k_A - k_0)$ for all finite $n$.*

If $\Sigma$ is orthogonally invariant, then given $\sigma$ and $\tau$, the LLR test statistic is the same as (18a) or (18b) with the norm defined by (12). As usual, when $\sigma^2$ and $\tau$ are unknown, they can be replaced by the consistent estimators of Lemma 3.3.

EXAMPLE 3.1 (Unrestricted). The whole space $\mathcal{M} = \mathcal{S}_p$ is trivially an affine subspace of itself and is parallel to $\mathcal{I}_p$. Here $\hat{M} = \bar{Y}$, $\hat{\sigma}^2 = s_\tau^2$ and

$$(19) \qquad \hat{\tau} = -\frac{\sum_{i=1}^n \|Y_i - \bar{Y}\|_{1, q/p}^2}{(q-1)\sum_{i=1}^n [\operatorname{tr}(Y_i - \bar{Y})]^2}.$$

Notice that the closer $\sum_{i=1}^n \operatorname{tr}[(Y_i - \bar{Y})^2]/q$ is to $\sum_{i=1}^n [\operatorname{tr}(Y_i - \bar{Y})]^2/p$, the closer is $\tau$ to 0, suggesting sphericity of the distribution. Here we have $\bar{Y} \sim N_{pp}(M, \sigma^2/n, \tau)$, and given $\tau$, $qns_\tau^2/\sigma^2 \sim \chi^2(q(n-1))$ independent of $\bar{Y}$.

The LLR test statistic for testing $H_0 : M = M_0$ vs. $H_A : M \neq M_0$ unrestricted is given immediately by Lemma 3.4. Since $\hat{M}_0 = M_0$ and $\hat{M}_A = \bar{Y}$, the LLR given $\sigma$ and $\tau$ is $T = n\|\bar{Y} - M_0\|_{\sigma^2, \tau}^2 \sim \chi^2(q)$ under $H_0$, and is independent of $\tau$. If $\sigma$ and $\tau$ are unknown, the LLR is an increasing function of the statistic $T = (n-1)\|\bar{Y} - M_0\|_{1, \hat{\tau}}^2 / (qs_\tau^2)$, which is approximately $F(q, q(n-1))$ under $H_0$ for large $n$. Both test statistics above reduce in the univariate case ($p = 1$, $\tau = 0$) to the squares of the known one-sample $z$ and $t$ statistics.

In addition to the above tests, one may want to test whether the orthogonally invariant covariance structure is an appropriate model for the data. This can be done using the LLR test given by the following proposition.

PROPOSITION 3.1. Let $Y_1, \ldots, Y_n \in \mathcal{S}_p$ i.i.d. $N_{pp}(M, \Sigma)$, $n > q(q+3)/2$, and consider the test $H_0 : \Sigma$ is orthogonally invariant with $\tau > 0$ vs. $H_A : \Sigma$ is unrestricted. The LLR statistic $2(l_A - l_0)$ for this test is

$$(20) \quad T = nq\log\hat{\sigma}^2 - n\log(1 - p\hat{\tau}) - n\log|\hat{\Sigma}| \xrightarrow[H_0]{n \to \infty} \chi^2\left(\frac{q(q+1)}{2} - 2\right),$$

where $\hat{\sigma}^2$ and $\hat{\tau}$ are given by (15) and (16), respectively, and $\hat{\Sigma}$ is the empirical covariance matrix $\hat{\Sigma} = (1/n)\sum_{i=1}^n \operatorname{vecd}(Y_i - \bar{Y})[\operatorname{vecd}(Y_i - \bar{Y})]'$.

PROOF. Replacing $\hat{\Sigma}$ in (14) gives the maximized likelihood under $H_A$, $l_A = -(n/2)\log|\hat{\Sigma}| - nq/2$ with dimension $q + q(q+1)/2$. Replacing $\hat{\sigma}^2$ and $\hat{\tau}$ in (17) gives the maximized likelihood under $H_0$, $l_0 = -(nq/2)\log\hat{\sigma}^2 + (n/2)\log(1 - p\hat{\tau}) - nq/2$ with dimension $q + 2$. The LLR $2(l_A - l_0)$ is equal to (20) and the difference of dimensions is $q(q+1)/2 - 2$. □

We emphasize that the asymptotic distribution in (20) is guaranteed only if $\tau < 0$ (i.e., $c > 0$). A more general situation of this kind is treated by Drton [8], Section 6. In that context, the set of orthogonally invariant covariances $\Sigma$ with $\tau \leq 0$ can be seen as a special case of a factor analysis model. The spherical case $\tau = 0$ is at the boundary of that semi-algebraic set and the asymptotic distribution of the LLR there is not $\chi^2$.

3.3. *General inference in the two-sample problem.* Let $Y_1, \ldots, Y_{n_1}$ and $Y_{n_1+1}, \ldots, Y_n$, $n = n_1 + n_2$, be two independent i.i.d. samples from $N_{pp}(M_1, \Sigma)$ and $N_{pp}(M_2, \Sigma)$, respectively. For simplicity, we assume the covariance $\Sigma$ is common to both samples, although this assumption can be relaxed in similar ways to the proposed solutions to the multivariate Behrens–Fisher problem (Scheffé [29], Mardia, Kent and Bibby [23]). Here we are interested in estimating and testing the pair $M = (M_1, M_2)$ when it is restricted to a subset $\mathcal{M}_2$ of $\mathcal{S}_p \times \mathcal{S}_p$ defined in terms of eigenvalues and eigenvectors of $M_1$ and $M_2$.

To solve two-sample problems it is convenient to define the following inner product in $\mathcal{S}_p \times \mathcal{S}_p$. For $A = (A_1, A_2)$ and $B = (B_1, B_2)$, let

$$\tag{21a} \langle\langle A, B \rangle\rangle_\Sigma = n_1 \langle A_1, B_1 \rangle_\Sigma + n_2 \langle A_2, B_2 \rangle_\Sigma$$

$$\tag{21b} = n \langle \mathrm{avg}(A), \mathrm{avg}(B) \rangle_\Sigma + \frac{n_1 n_2}{n} \langle \Delta(A), \Delta(B) \rangle_\Sigma,$$

where we have defined for a generic $X = (X_1, X_2)$ the linear functions

$$\mathrm{avg}(X) = \frac{n_1 X_1 + n_2 X_2}{n}, \qquad \Delta(X) = X_1 - X_2.$$

Define the group averages $\bar{Y}_1 = (1/n_1) \sum_{i=1}^{n_1} Y_i$ and $\bar{Y}_2 = (1/n_2) \sum_{i=n_1+1}^{n} Y_i$, and let $\bar{Y} = (\bar{Y}_1, \bar{Y}_2)$. Writing the joint likelihood of the sample and pooling terms leads to the MLE of $M = (M_1, M_2)$ for general sets $\mathcal{M}_2 \subset \mathcal{S}_p \times \mathcal{S}_p$ as given by the following lemma.

LEMMA 3.5. *Suppose the pair $M = (M_1, M_2)$ lies in a nonempty closed subset $\mathcal{M}_2$ of $\mathcal{S}_p \times \mathcal{S}_p$. Then, given $\Sigma$, the MLE $\hat{M} = (\hat{M}_1, \hat{M}_2)$ of $M$ is unique almost surely and minimizes over $\mathcal{M}_2$ the squared distance*

$$\tag{22a} g_\Sigma(M_1, M_2) = n_1 \|\bar{Y}_1 - M_1\|_\Sigma^2 + n_2 \|\bar{Y}_2 - M_2\|_\Sigma^2$$

$$\tag{22b} = n\|\mathrm{avg}(\bar{Y}) - \mathrm{avg}(M)\|_\Sigma^2 + \frac{n_1 n_2}{n} \|\Delta(\bar{Y}) - \Delta(M)\|_\Sigma^2.$$

*The MLE $\hat{M}$ is consistent as $n_1, n_2 \to \infty$.*

As in the one-sample case, we show that for a class of sets $\mathcal{M}_2 \subset \mathcal{S}_p \times \mathcal{S}_p$, the MLE does not depend on $\Sigma$ if and only if $\Sigma$ is orthogonally invariant.

DEFINITION 3.2. Let $\mathcal{I}_p \times \mathcal{I}_p = \{(M_1, M_2) = (\alpha I_p, \beta I_p), \alpha, \beta \in \mathbb{R}\}$ and $(\mathcal{I}_p, \mathcal{I}_p) = \{(M_1, M_2) = \gamma(I_p, I_p), \gamma \in \mathbb{R}\} \subset \mathcal{I}_p \times \mathcal{I}_p$. We say that a set $\mathcal{M}_2 \subset \mathcal{S}_p \times \mathcal{S}_p$ is:

(a) orthogonal to $\mathcal{I}_p \times \mathcal{I}_p$, if $\mathcal{M}_2$ can be written as $\mathcal{M}_2 = \mathcal{A} \oplus \{(aI_p, bI_p)\}$ for some $a, b \in \mathbb{R}$ and $\mathcal{A} \subset \mathcal{S}_p \times \mathcal{S}_p$ such that $\mathrm{tr}(A_1) = \mathrm{tr}(A_2) = 0$ for all $(A_1, A_2) \in \mathcal{A}$.

(b) parallel to $(\mathcal{I}_p, \mathcal{I}_p)$, if it can be written as $\mathcal{M}_2 = \mathcal{A} \oplus (\mathcal{I}_p, \mathcal{I}_p)$, where $\mathcal{A}$ is orthogonal to $\mathcal{I}_p \times \mathcal{I}_p$ in the sense of definition (a).

(c) parallel to $\mathcal{I}_p \times \mathcal{I}_p$, if it can be written as $\mathcal{M}_2 = \mathcal{A} \oplus (\mathcal{I}_p \times \mathcal{I}_p)$, where $\mathcal{A}$ is orthogonal to $\mathcal{I}_p \times \mathcal{I}_p$ in the sense of definition (a).

The definition of orthogonality Definition 3.2(a) is equivalent to the condition that for every tangent vector $A = (A_1, A_2)$ to $\mathcal{M}_2$, $\langle\!\langle A, (\alpha I_p, \beta I_p) \rangle\!\rangle_{\sigma^2, \tau} = 0$ for all $\alpha, \beta \in \mathbb{R}$. This can be checked easily replacing in (21b): setting $\alpha = \beta$ gives $\operatorname{tr}(\operatorname{avg}(A)) = 0$ and using $\alpha = -n_2/n_1 \beta$ gives $\operatorname{tr}(\Delta(A)) = 0$. All three properties defined in Definition 3.2 are orthogonally invariant: if $\mathcal{M}_2$ belongs to one of the three categories, the set $\{(Q_1 M_1 Q_1', Q_2 M_2 Q_2'), (M_1, M_2) \in \mathcal{M}_2\}$, where $Q_1, Q_2 \in \mathcal{O}_p$, also belongs to the same category as $\mathcal{M}_2$. It is easy to check that each of the sets in Figure 3 belongs to one of the categories in Definition 3.2(b) or 3.2(c): the affine subset $\{M_1 = M_2\}$ and the set $M_{2,D}$ are parallel to $(\mathcal{I}_p, \mathcal{I}_p)$; the other two are parallel to $\mathcal{I}_p \times \mathcal{I}_p$.

THEOREM 3.2.  $\Sigma$ *is orthogonally invariant if and only if, for every embedded submanifold* $\mathcal{M}_2 \subset \mathcal{S}_p \times \mathcal{S}_p$ *that is either parallel to* $(\mathcal{I}_p, \mathcal{I}_p)$ *or parallel to* $\mathcal{I}_p \times \mathcal{I}_p$, *every critical point* $\hat{M} = (\hat{M}_1, \hat{M}_2)$ *of* $g_\Sigma(M) = n_1 \|\bar{Y}_1 - M_1\|_\Sigma^2 + n_2 \|\bar{Y}_2 - M_2\|_\Sigma^2$ *over* $\mathcal{M}_2$ *is also a critical point of* $g_{1,0}(M) = n_1 \operatorname{tr}[(\bar{Y}_1 - M_1)^2] + n_2 \operatorname{tr}[(\bar{Y}_2 - M_2)^2]$.

PROOF.  (i) Suppose $\Sigma$ is orthogonally invariant and let $B = \bar{Y} - \hat{M} \in \mathcal{S}_p \times \mathcal{S}_p$. Every critical point $\hat{M}$ of $g_\Sigma(M)$ must satisfy the orthogonality condition $\langle\!\langle A, B \rangle\!\rangle_{\sigma^2, \tau} = 0$ for every tangent vector $A = (A_1, A_2)$ to $\mathcal{M}_2$ at $\hat{M}$. Using (12) and (21b), the orthogonality condition can be written as

$$
\begin{aligned}
(23) \qquad & \left[ n \operatorname{tr}[\operatorname{avg}(A) \operatorname{avg}(B)] + \frac{n_1 n_2}{n} \operatorname{tr}[\Delta(A)\Delta(B)] \right] \\
& - \tau \left[ n \operatorname{tr}[\operatorname{avg}(A)] \operatorname{tr}[\operatorname{avg}(B)] + \frac{n_1 n_2}{n} \operatorname{tr}[\Delta(A)] \operatorname{tr}[\Delta(B)] \right] = 0.
\end{aligned}
$$

$\hat{M}$ is a critical point of $g_{1,0}(M)$ only if the first bracket in (23) is zero. We thus need to show that the second bracket is equal to zero for all $\tau$ and all tangent vectors $A$.

Assume $\mathcal{M}_2$ is parallel to $(\mathcal{I}_p, \mathcal{I}_p)$. The tangent space to $\mathcal{M}_2$ at $\hat{M}$ can be written as $\mathcal{A} \oplus (\mathcal{I}_p, \mathcal{I}_p)$ where every tangent vector $A = (A_1, A_2) \in \mathcal{A}$ satisfies $\operatorname{tr}(A_1) = \operatorname{tr}(A_2) = 0$. For every $A \in \mathcal{A}$, $\operatorname{tr}(\operatorname{avg}(A)) = \operatorname{tr}(\Delta(A)) = 0$ and we are done. For the remaining component parallel to $(\mathcal{I}_p, \mathcal{I}_p)$, let $A = \alpha(I_p, I_p)$ so that $\operatorname{avg}(A) = \alpha I_p$ and $\Delta(A) = 0$. Replacing in (23), we see that $\operatorname{tr}[\operatorname{avg}(B)] = 0$ for any critical point with respect to $\langle\!\langle \cdot, \cdot \rangle\!\rangle_{\sigma^2, \tau}$. Since $\Delta(A) = 0$, the second bracket in (23) vanishes for all $\tau$.

Assume instead $\mathcal{M}_2$ is parallel to $\mathcal{I}_p \times \mathcal{I}_p$. The tangent space to $\mathcal{M}_2$ at $\hat{M}$ can be written as $\mathcal{A} \oplus (\mathcal{I}_p, \mathcal{I}_p) \oplus (n_2 \mathcal{I}_p, -n_1 \mathcal{I}_p)$ where $\mathcal{A}$ is orthogonal to $\mathcal{I}_p \times \mathcal{I}_p$ and $(n_2 \mathcal{I}_p, -n_1 \mathcal{I}_p) = \{(M_1, M_2) = \gamma(n_2 I_p, -n_1 I_p), \gamma \in \mathbb{R}\} \subset \mathcal{I}_p \times \mathcal{I}_p$. The first two components were covered in the previous case. For every tangent vector $A$ in the component parallel to $(n_2 \mathcal{I}_p, -n_1 \mathcal{I}_p)$, let $A = \alpha(n_2 I_p, -n_1 I_p)$, so that $\mathrm{avg}(A) = 0$ and $\Delta(A) = \alpha n I_p$. Replacing in (23), we see that $\mathrm{tr}[\Delta(B)] = 0$ for any critical point with respect to $\langle\!\langle \cdot, \cdot \rangle\!\rangle_\Sigma$. Since $\mathrm{avg}(A) = 0$, the second bracket in (23) vanishes for all $\tau$.

(ii) The proof is similar to the proof in the one-sample case (Theorem 3.1). $\square$

By Theorem 3.2, if we assume the orthogonally invariant model (11), then within the class of sets $\mathcal{M}_2$ that are either parallel to $(\mathcal{I}_p, \mathcal{I}_p)$ or parallel to $\mathcal{I}_p \times \mathcal{I}_p$, the problem of finding the MLE of $M$ reduces to minimizing $g_{1,0}(M)$ over $\mathcal{M}_2$. Once $\hat{M}$ is found, the estimates of $\sigma$ and $\tau$ in model (11) easily follow as shown in the following lemma.

LEMMA 3.6. *Let* $Y_1, \ldots, Y_{n_1}$ *and* $Y_{n_1+1}, \ldots, Y_n$, $n = n_1 + n_2$, *be two independent i.i.d. samples from* $N_{pp}(M_1, \sigma^2, \tau)$ *and* $N_{pp}(M_2, \sigma^2, \tau)$, *respectively,* $p \geq 2$. *Suppose the pair* $M = (M_1, M_2)$ *lies in a nonempty closed subset* $\mathcal{M}_2$ *of* $\mathcal{S}_p \times \mathcal{S}_p$ *and denote by* $\hat{M} = (\hat{M}_1, \hat{M}_2)$ *the MLE of* $M$. *Then, given* $\hat{M}$ *and* $\tau$, *the MLE of* $\sigma^2$ *is given by*

$$(24) \qquad \hat{\sigma}^2 = s_{12}^2 + \frac{1}{qn}(n_1 \|\bar{Y}_1 - \hat{M}_1\|_{\sigma^2, \tau}^2 + n_2 \|\bar{Y}_2 - \hat{M}_2\|_{\sigma^2, \tau}^2),$$

*where*

$$s_{12}^2 = \frac{1}{qn}\left(\sum_{i=1}^{n_1} \|Y_i - \bar{Y}_1\|_{1,\tau}^2 + \sum_{i=n_1+1}^{n} \|Y_i - \bar{Y}_2\|_{1,\tau}^2\right)$$

*is the pooled variance. Moreover, given* $\hat{M}$, *the MLE of* $\tau$ *is*

$$\hat{\tau} = -\frac{\sum_{i=1}^{n} \|Y_i - \mathrm{avg}(\bar{Y})\|_{1,q/p}^2 + n_1 \|\bar{Y}_1 - \hat{M}_1\|_{1,q/p}^2 + n_2 \|\bar{Y}_2 - \hat{M}_2\|_{1,q/p}^2}{(q-1)\{\sum_{i=1}^{n}[\mathrm{tr}(Y_i - \mathrm{avg}(\bar{Y}))]^2 + n_1[\mathrm{tr}(\bar{Y}_1 - \hat{M}_1)^2] + n_1[\mathrm{tr}(\bar{Y}_2 - \hat{M}_2)^2]\}}.$$

*The above estimates are unique and are consistent as* $n_1, n_2 \to \infty$.

In the testing situation, let $Y_1, \ldots, Y_{n_1}$ and $Y_{n_1+1}, \ldots, Y_n$, $n = n_1 + n_2$, be two independent i.i.d. samples from $N_{pp}(M_1, \Sigma)$ and $N_{pp}(M_2, \Sigma)$, respectively. Consider the test $H_0 : (M_1, M_2) \in \mathcal{M}_0$ vs. $H_0 : (M_1, M_2) \in \mathcal{M}_A$, where $\mathcal{M}_0 \subset \mathcal{M}_A \subseteq \mathcal{S}_p \times \mathcal{S}_p$ are closed embedded submanifolds with dimensions $k_A = \dim(\mathcal{M}_A)$, $k_0 = \dim(\mathcal{M}_0)$, $k_A > k_0$. Let $\hat{M}_{1,0}$, $\hat{M}_{2,0}$ denote the MLEs of $M_1$ and $M_2$ under $H_0$, and $\hat{M}_{1,A}$, $\hat{M}_{2,A}$ the corresponding MLEs under $H_A$.

Lemma 3.7. *Suppose $\mathcal{M}_0 \subset \mathcal{M}_A$ are close embedded submanifolds of $\mathcal{S}_p \times \mathcal{S}_p$. Given $\Sigma$, the LLR test statistic $2(l_A - l_0)$ is*

$$
\begin{aligned}
(25) \qquad T = {}& n_1 \|\bar{Y}_1 - \hat{M}_{1,0}\|_\Sigma^2 + n_2 \|\bar{Y}_2 - \hat{M}_{2,0}\|_\Sigma^2 \\
& - n_1 \|\bar{Y}_1 - \hat{M}_{1,A}\|_\Sigma^2 - n_2 \|\bar{Y}_2 - \hat{M}_{2,A}\|_\Sigma^2
\end{aligned}
$$

*and is distributed as $\chi^2(k_A - k_0)$ under $\mathcal{M}_0$ asymptotically as $n_1, n_2 \to \infty$. If $\mathcal{M}_0 \subset \mathcal{M}_A$ are both affine subspaces of $\mathcal{S}_p \times \mathcal{S}_p$, then the distribution is exactly $\chi^2(k_A - k_0)$ for all finite $n_1$ and $n_2$.*

If $\Sigma$ is orthogonally invariant, then given $\sigma$ and $\tau$, the LLR test statistic $2(l_A - l_0)$ is the same as (25) with the norm defined by (12). As in the one-sample case, if $\sigma^2$ and $\tau$ are unknown, then they can be replaced by the consistent estimators of Lemma 3.6.

Example 3.2 (Unrestricted). If $\mathcal{M}_2 = \mathcal{S}_p \times \mathcal{S}_p$ unrestricted, the MLEs are $\hat{M}_1 = \bar{Y}_1$, $\hat{M}_2 = \bar{Y}_2$, $\hat{\sigma}^2 = s_{12}^2$ and $\hat{\tau}$ is given by (19). $\bar{Y}_1$ and $\bar{Y}_2$ are independent $N_{pp}(M_1, \sigma^2/n, \tau)$ and $N_{pp}(M_2, \sigma^2/n, \tau)$, respectively, and are independent of $s_{12}^2$ given $\tau$. In addition, given $\tau$, $qns_{12}^2/\sigma^2 \sim \chi^2(q(n-2))$.

The LLR test statistic for testing $H_0 : M_1 = M_2$ vs. $H_A : M_1 \neq M_2$ is given by Lemma 3.7. Replacing $\hat{M} = \hat{M}_{1,0} = \hat{M}_{2,0} = \bar{Y}$, $\hat{M}_{1,A} = \bar{Y}_1$, and $\hat{M}_{2,A} = \bar{Y}_2$ in (25), the test statistic given $\sigma$ and $\tau$ is $T = (n_1 n_2/n)\|\bar{Y}_1 - \bar{Y}_2\|_{\sigma^2, \tau}^2 \sim \chi^2(q)$ under $H_0$. If $\sigma$ and $\tau$ are unknown, then the LLR is an increasing function of the statistic $T = (n-2)n_1 n_2 \|\bar{Y}_1 - \bar{Y}_2\|_{1,\hat{\tau}}^2/(qn^2 s_{12}^2)$ which is approximately $F(q, q(n-2))$ under $H_0$ for large $n_1$ and $n_2$. Both test statistics above reduce in the univariate case ($p = 1$, $\tau = 0$) to the squares of the known two-sample $z$ and $t$-statistics.

## 4. Inference for eigenvalues and eigenvectors in the one-sample problem.
For any of the sets described in Section 2.1 and Figure 1, $\hat{M}$ immediately specifies the MLE of $\sigma^2$ and $\tau$ by Lemma 3.3. In addition, by Lemma 3.4, all we need to know is the form of the test statistic (18a) or (18b) and corresponding test statistics for unknown $\sigma$ and $\tau$ may be derived from there. Thus, for brevity, we assume from here on that $\sigma$ and $\tau$ are known. To fix ideas, we begin with the linear cases and then move on to the nonlinear ones.

### 4.1. *Affine subspaces.*
The unrestricted case (A0) was covered in Example 3.1. Both cases (A1) and (A2) involve the same affine subspace $\mathcal{M}_{U_0}$ defined by (2). In order to derive the LLR statistics for both these cases we need the MLE of $M$ in $\mathcal{M}_{U_0}$.

PROPOSITION 4.1. *Let $\mathcal{M}_{U_0}$ be given by* (2). *The MLE of $M$ over $\mathcal{M}_{U_0}$ is $\hat{M} = U_0 \hat{D} U_0'$ where*

$$\hat{D} = \mathrm{diag}(U_0' \bar{Y} U_0). \tag{26}$$

*The diagonal elements of $\hat{D}$ are multivariate normal $N(\mathrm{diag}(D), (I_p + c\mathbf{1}_p\mathbf{1}_p')\sigma^2/n)$, where $c = \tau/(1 - p\tau)$.*

PROOF. Because $\mathcal{M}_{U_0}$ is an affine subspace, the MLE of $M$ is found by orthogonal projection. $M_{U_0}$ is parallel to $\mathcal{I}_p$. By Theorem 3.1, $\hat{M} = U_0 \hat{D} U_0'$, where $\hat{D}$ is the minimizer of $g_{1,0}(D) = \mathrm{tr}[(\bar{Y} - U_0 D U_0')^2] = \mathrm{tr}[(U_0' \bar{Y} U_0 - D)^2]$. The solution is the diagonal matrix $\hat{D}$ closest in Frobenius norm to $U_0' \bar{Y} U_0$, that is, (26). It is easy to see that $\hat{D}$ is unbiased. Also, because of the orthogonal invariance property, $U_0' \bar{Y} U_0 \sim N_{pp}(U_0' M U_0, \sigma^2/n, \tau) = N_{pp}(D, \sigma^2/n, \tau)$, so $\mathrm{diag}(\hat{D}) \sim N(\mathrm{diag}(D), (I_p + c\mathbf{1}_p\mathbf{1}_p')\sigma^2/n)$. □

The estimate $\hat{D}$ is invariant under sign changes of the eigenvectors and depends on the order of the columns of $U_0$ only in the sense that permuting the columns of $U_0$ will result in the same permutation of the columns of $\hat{D}$. If $P$ is a signed permutation matrix with elements $+1$ or $-1$ and the eigenvector matrix $U_0$ is permuted to $U_0 P$, then $\mathrm{diag}((U_0 P)' \bar{Y} (U_0 P)) = \mathrm{diag}(P'(U_0' \bar{Y} U_0)P) = P' \mathrm{diag}(U_0' \bar{Y} U_0)P = P' \hat{D} P$. Notice that $\hat{D}$ is the same whether or not the true underlying eigenvalues $D$ have multiplicities.

In case (A1), the LLR test statistic for testing $H_0 : M = M_0$ vs. $H_A : M \in \mathcal{M}_{U_0}$ is obtained from Lemma 3.4. Here $\hat{M}_0 = U_0 D_0 U_0'$ and $\hat{M}_A = U_0 \hat{D} U_0'$, where $\hat{D} = \mathrm{diag}(U_0' \bar{Y} U_0)$, so given $\sigma$ and $\tau$, the LLR test statistic is

$$T = n\|\hat{D} - D_0\|_{\sigma^2, \tau}^2 \underset{H_0}{\sim} \chi^2(p). \tag{27}$$

In case (A2), the test is $H_0 : M \in \mathcal{M}_{U_0}$ vs. $M \notin \mathcal{M}_{U_0}$. Here $\hat{M}_0 = U_0 \hat{D} U_0'$, $\hat{D} = \mathrm{diag}(U_0' \bar{Y} U_0)$, and $\hat{M}_A = \bar{Y}$, so given $\sigma$ and $\tau$, the LLR test statistic is

$$T = n\|\bar{Y} - U_0 \hat{D} U_0'\|_{\sigma^2, \tau}^2 \underset{H_0}{\sim} \chi^2(q - p). \tag{28}$$

4.2. *A convex polyhedral cone [Case (C2)].* Let $\mathcal{M}_{U_0}^{>}$ be defined by (3). This set is parallel to $\mathcal{I}_p$ according to Definition 3.1, and Theorem 3.1 applies by checking each face for critical points. Define the vectors $y = \mathrm{diag}(U_0' \bar{Y} U_0)$ and $d = \mathrm{diag}(D)$. A similar argument as in the proof of Proposition 4.1 gives that the MLE of $M$ is $\hat{M} = U_0 \hat{D} U_0'$, where the vector of diagonal elements $\hat{d}$ of $\hat{D}$ is the minimizer of $|y - \hat{d}|^2$ constrained to $\hat{d}_1 \geq \cdots \geq \hat{d}_p$. Many algorithms have been proposed for solving this type of problem, sometimes called order-restricted least squares, isotonic least squares or nonnegative least squares (since the consecutive differences are nonnegative) (Lawson and Hanson [19], Dykstra [9]). A direct way to solve it in our case is to apply the pool-adjacent-violators algorithm (PAVA) with equal initial weights (Robertson and Wegman [28], Raubertas [27]).

If the true parameter $D$ is such that $d_1 > \cdots > d_p$, then $\sqrt{n}(\hat{d} - d)$ is asymptotically multivariate normal $N(0, \sigma^2(I_p + c\mathbf{1}_p\mathbf{1}_p'))$. Otherwise, if $d$ has multiplicities, the asymptotic distribution of $\sqrt{n}(\hat{d} - d)$ is the projection of $N(0, \sigma^2(I_p + c\mathbf{1}_p\mathbf{1}_p'))$ onto the support cone at the true $d$ (Self and Liang [31]). For example, if $p = 3$ and the true parameter $d$ is such that $d_1 = d_2 > d_3$, then the support cone at $d$ is the cone $d_1 \geq d_2$ (because the inequality $d_2 > d_3$ asymptotically has no effect). The asymptotic distribution of $\sqrt{n}(\hat{d} - d)$ is the projection of $N(0, \sigma^2(I_3 + c\mathbf{1}_3\mathbf{1}_3'))$ onto the support cone at $d$: it has mass $1/2$ on the plane $d_1 = d_2$ and is normal on the half-space $d_1 > d_2$.

PROPOSITION 4.2 [Case (C2)].  *Consider the test* $H_0: M \in \mathcal{M}_{U_0}^{>}$ *vs.* $M \notin \mathcal{M}_{U_0}^{>}$, *and let* $U_0\hat{D}U_0'$ *be the MLE of* $M$ *in* $\mathcal{M}_{U_0}^{>}$. *Given* $\sigma$ *and* $\tau$, *the LLR statistic is* $T = n\|\bar{Y} - U_0\hat{D}U_0'\|_{\sigma^2, \tau}^2$. *If the true* $D$ *has* $k \leq p$ *distinct eigenvalues, then the asymptotic distribution of* $T$ *is a mixture of* $p - k + 1$ $\chi^2$-*variables with degrees of freedom* $q - p, q - p + 1, \ldots, q - k$. *The weights of the mixture depend only on the angles of the cone at* $M = U_0DU_0'$ *and do not depend on* $\sigma$ *or* $\tau$.

PROOF.  The expression for $T$ is obtained in a similar way to (28). Its asymptotic distribution follows from the theory of LLRs on convex cones (Self and Liang [31]). The weights of the mixture do not depend on $\sigma$ because $\sigma$ is a scaling factor that cancels out in $T$. They also do not depend on $\tau$ because $\mathcal{M}_{U_0}^{>}$ is parallel to $\mathcal{I}_p$ and the relevant angles between the faces of the cone are computed as inner products on the space orthogonal to $\mathcal{I}_p$. These inner products do not depend on $\tau$ by Lemma 3.1.  □

As specific cases, if $D$ has $p$ distinct eigenvalues, then $T$ is asymptotically $\chi^2(q - p)$. The largest number of mixture components possible is $p$, obtained when $D$ is isotropic. In this case we get a mixture of $\chi^2$ variables with degrees of freedom $q - p, \ldots, q - 1$. As an example, take $p = 3$. If the true underlying $d$ is oblate ($d_1 = d_2 > d_3$), then the cone at $d$ locally looks like a half space and the asymptotic distribution of $T$ is $(1/2)\chi_3^2 + (1/2)\chi_4^2$. On the other hand, if $D$ is isotropic ($d_1 = d_2 = d_3$), the cone at $d$ is framed by the intersection of two planes at an angle of $60°$, so the asymptotic distribution of $T$ is $(1/6)\chi_3^2 + (1/2)\chi_4^2 + (1/3)\chi_5^2$.

For general $p$, the mixture weights can be easily obtained by simulation. Given $d$, generate $y \sim N(d, I_p)$ and project it onto the appropriate cone using the PAVA algorithm while keeping track of the dimension of the face that $y$ is projected onto. Repeating this many times, the weight corresponding to the $\chi^2(q - k)$ mixture component is the proportion of times that $y$ gets projected onto a face of dimension $k$.

4.3. *Curved submanifolds [Cases (S1) and (S2)].* Both tests (S1) and (S2) involve the same submanifold $\mathcal{M}_{D_0}$. The dimension of the submanifold $\mathcal{M}_{D_0}$ depends on the multiplicities of $D_0$ and is given by the following lemma. The MLE of $M$ in $\mathcal{M}_{D_0}$ is then given by Theorem 4.1 below.

LEMMA 4.1.   $\dim(\mathcal{M}_{D_0}) = q - \sum_{i=1}^{k} m_i(m_i + 1)/2$.

PROOF.   $\mathcal{M}_{D_0}$ is diffeomorphic to the quotient $\mathcal{O}_p / O m_1 \times \cdots \times O m_k$. This is because if $D_0$ has diagonal blocks $d_1 I_{m_1}, \ldots, d_k I_{m_k}$ and $Q$ is block diagonal orthogonal with orthogonal diagonal blocks $Q_1, \ldots, Q_k$, then $Q' D_0 Q = D_0$. The set of such matrices $Q$ is $O m_1 \times \cdots \times O m_k$. Therefore,

$$\dim(\mathcal{M}) = \dim(\mathcal{O}_p) - \sum_{i=1}^{k} \dim(O(m_i)) = \frac{p(p-1)}{2} - \frac{1}{2} \sum_{i=1}^{k} m_i(m_i - 1). \quad \square$$

THEOREM 4.1.   *Let $\mathcal{M}_{D_0}$ be given by (4) and let $\bar{Y} = V \Lambda V'$ be an eigendecomposition chosen so that the diagonal elements of $\Lambda$ are in decreasing order.*

(i) *The MLE of $M$ is $\hat{M} = \hat{U} D_0 \hat{U}'$, where $\hat{U}$ is any matrix of the form*

(29) $$\hat{U} = VQ,$$

*and $Q$ is an orthogonal matrix such that $Q D_0 Q' = D_0$, that is, block diagonal orthogonal with diagonal orthogonal blocks of sizes $m_1, \ldots, m_k$. In particular, if the eigenvalues of $D_0$ are distinct, then such $Q$'s are diagonal matrices with diagonal entries equal to $\pm 1$.*

(ii) *Assume $D_0$ has distinct diagonal entries $d_1 > \cdots > d_p$, in which case $M$ has unique eigenvectors $U$ up to sign. Let $\hat{U}$ in (29) be chosen to minimize $\|U \hat{U}' - I_p\|$ and let $\hat{A} = \log(U' \hat{U}) \in \mathcal{A}_p$, the set of $p \times p$ antisymmetric matrices. Then, as $n \to \infty$, the entries $\{\hat{a}_{ij}\}_{i \neq j}$ of $\hat{A}$ are asymptotically independent and*

$$\sqrt{n} \hat{a}_{ij} \to N\left(0, \frac{\sigma^2}{2(d_i - d_j)^2}\right).$$

Computing $\hat{U}$ above leads to a maximization problem in $\mathcal{O}_p$. We solve it by requiring that the gradient of the objective function be orthogonal to the tangent space to $\mathcal{O}_p$. For this, we use the fact that the tangent vectors $\dot{U}$ at any point $U \in \mathcal{O}_p$ are of the form $\dot{U} = UA$ where $A \in \mathcal{A}_p$, the set of $p \times p$ antisymmetric matrices (Chang [5], Edelman, Arias and Smith [10], Lang [18], Moakher [25]). Part (ii) of the theorem gives the asymptotic distribution of the error incurred in the estimation of $U$. However, instead of using the Euclidean difference $\hat{U} - U$ to measure the error, we use the Riemannian logarithmic map of $\hat{U}$ to the tangent space to $\mathcal{O}_p$ at $U$, which results in the tangent vector $\hat{A} = \log(U'\hat{U}) \in \mathcal{A}_p$ (Edelman, Arias and Smith [10], Lang [18], Moakher [25]). Notice that the variance of $\hat{a}_{ij}$ increases without bound as $d_i$ and $d_j$ get closer. This happens because the curvature of $\mathcal{M}_{D_0}$ on the plane $d_i + d_j = $ constant increases without bound as $d_i$ and $d_j$ get closer. In the limit when $d_i = d_j$, we get one eigenvalue with multiplicity two and the dimension of the set abruptly goes down by 1. In Figure 2, this phenomenon can be seen as the collapse of the circle into its center. The result is that estimates of eigenvectors become increasingly variable as the corresponding eigenvalues get close, and unidentifiable if they are equal.

PROOF OF THEOREM 4.1. (i) $\mathcal{M}_{D_0}$ is parallel to $\mathcal{I}_p$. Thus by Theorem 3.1, $\hat{U}$ is the minimizer of

$$(30) \qquad g(U) = \operatorname{tr}[(\bar{Y} - UD_0U')^2] = \operatorname{tr}(\bar{Y}^2) - 2\operatorname{tr}(U'\bar{Y}UD_0) + \operatorname{tr}(D_0^2)$$

subject to $UU' = I_p$. The critical points of $g(U)$ with respect to the constraint $U'U = I_p$ are those points $\hat{U}$ where the gradient is orthogonal to the surface $U'U = I_p$. The tangent vectors to that surface at $U$ are of the form $\dot{U} = UA$, where $A \in \mathcal{A}_p$. To see this, trace a curve $Q(t) \in \mathcal{O}(p)$ such that $Q(0) = U$. Differentiating $U'U = I_p$ with respect to $t$ and evaluating at $t = 0$ gives $\dot{U}'U + U'\dot{U} = 0$, so $A = U'\dot{U} \in \mathcal{A}_p$.

Using matrix derivative rules (Fang and Zhang [12], page 16) we get that the gradient of $g(U)$ is $\partial g(U)/\partial U = -4\bar{Y}UD_0$. The gradient must satisfy

$$\left\langle \dot{U}, \frac{\partial g}{\partial U} \right\rangle_{U=\hat{U}} = -4\operatorname{tr}[(\hat{U}A)'\bar{Y}\hat{U}D_0] = -4\operatorname{tr}(A'\hat{U}'\bar{Y}\hat{U}D_0) = 0$$

for all $A \in \mathcal{A}_p$. It is easy to check that the space $\mathcal{A}_p$ is orthogonal to $\mathcal{S}_p$, that is, $\operatorname{tr}(AB) = 0$ for all $A \in \mathcal{A}_p$ and $B \in \mathcal{S}_p$. Thus we must have that $\hat{U}'\bar{Y}\hat{U}D_0$ is symmetric. This is satisfied by $\hat{U} = VQ$, where $Q$ is any orthogonal matrix such that $QD_0 = D_0Q$. Plugging back into (30) gives $g(\hat{U}) = \operatorname{tr}[(\bar{Y} - VQD_0Q'V')^2] = \operatorname{tr}[(\Lambda - D_0)^2]$, which is minimized if the eigenvalues in $\Lambda$ are chosen in decreasing order. Of course, eigenvalues in $\Lambda$ corresponding to sets of repeated eigenvalues in $D_0$ can be permuted without affecting the minimization.

(ii) Assume momentarily that the covariance structure is spherical ($\tau = 0$). The definition $\hat{A} = \log(U'\hat{U})$ provides a parametrization of $\mathcal{SO}_p$ (the set of orthogonal matrices of determinant 1) given by $\hat{U} = U\exp(\hat{A})$, while at the same time $\hat{A}$ parametrizes the tangent space to $\mathcal{SO}_p$ at $U$. Minus twice the maximized log-likelihood is $ng(\hat{U})$. To get the covariance of $\hat{A}$, we need the score function, which is the gradient of $ng(\hat{U})$. To obtain it, trace a curve $Q(t) = U\exp(\hat{A}t)$ passing through $U$ in the direction $\hat{A}$. Let $W = U'\bar{Y}U$. The derivative of the log-likelihood in the direction $\hat{A}$ is the derivative of $ng(Q(t))$ with respect to $t$ evaluated at $t = 0$:

$$
\begin{aligned}
n\frac{dg(Q(t))}{dt}\Big|_{t=0} &= -2n\frac{d}{dt}\operatorname{tr}(e^{\hat{A}'t}We^{\hat{A}t}D)|_{t=0} \\
&= -2n\operatorname{tr}(e^{\hat{A}'t}\hat{A}'We^{\hat{A}t}D + e^{\hat{A}'t}W\hat{A}e^{\hat{A}t}D)|_{t=0} \\
&= -2n\operatorname{tr}(\hat{A}'WD + W\hat{A}D) = 2n\operatorname{tr}[\hat{A}'(DW - WD)].
\end{aligned}
$$

This is the inner product of $\hat{A}$ with the antisymmetric matrix $S = 2n(DW - WD)$. Thus $S$ is the score function in the standard basis for $\mathcal{A}_p$ and $\mathcal{S}_p$. It is easy to see that the entries of $S$ are $s_{ij} = 2n(d_i - d_j)w_{ij}$ for $i \neq j$. Since $W = U'\bar{Y}U \sim N_{pp}(D, \sigma^2/n)$, the off-diagonal entries of $W$ are independent $N(0, \sigma^2/(2n))$. Thus the $s_{ij}$ are independent normal with mean 0 and variance $\operatorname{var}(s_{ij}) = 4n(d_i - d_j)^2\operatorname{var}(w_{ij}) = 2n(d_i - d_j)^2\sigma^2$. This is the Fisher information with respect to $a_{ij}$. Thus, asymptotically, the entries $\sqrt{n}\hat{a}_{ij}$ are independent and each normal with mean zero and variance $\sigma^2/(2(d_i - d_j)^2)$.

For the above calculations we assumed $\tau = 0$. However, the tangent space to $\mathcal{M}_{D_0}$ is orthogonal to $\mathcal{I}_p$ according to Definition 3.1. By Lemma 3.1, all the inner products computed above are the same if $\tau \neq 0$. Therefore the asymptotic distribution obtained above is also the same if $\tau \neq 0$. $\quad\square$

We can now compute the LLR statistics for tests (S1) and (S2).

COROLLARY 4.1 [Case (S1)]. *Consider the test* $H_0: M = M_0$ *vs.* $M \in \mathcal{M}_{D_0}$. *Let* $\bar{Y} = V\Lambda V'$ *be an eigendecomposition of* $\bar{Y}$ *so that* $\Lambda$ *and* $D_0$ *are in decreasing order. Given* $\sigma$ *and* $\tau$, *the LLR test statistic is*

$$
(31) \quad T = \frac{2n}{\sigma^2}[\operatorname{tr}(\Lambda D_0) - \operatorname{tr}(\bar{Y}M_0)] \xrightarrow[H_0]{n\to\infty} \chi^2\left(q - \frac{1}{2}\sum_{i=1}^{k}m_i(m_i + 1)\right),
$$

*independent of* $\tau$.

PROOF. Here $\hat{M}_0 = M_0$ and $\hat{M}_A = VD_0V'$. Replacing these in (18b) gives $T = 2n\langle\bar{Y}, VD_0V' - M_0\rangle_{\sigma^2,\tau}$. But the second argument of the inner product has zero trace so the $\tau$-term drops. Thus $T = 2n\operatorname{tr}[\bar{Y}(VD_0V' - M_0)]/\sigma^2$, which equals (31). The number of degrees of freedom is equal to the dimension of $\mathcal{M}_{D_0}$ minus zero for the single point $M_0$. $\quad\square$

The form of the test statistic (31) is interesting. Recall that (S1) is a test of whether $M_0$ has eigenvectors $U_0$ when the eigenvalues $D_0$ are assumed known. If the eigenvectors $V$ of $\bar{Y}$ are equal (up to sign) to the eigenvectors $U_0$ of $M_0$, then $T$ is equal to zero. As the angles between the eigenvectors in $V$ and $U_0$ increase, the inner product between $\bar{Y}$ and $M_0$ decreases, which increases the value of the test statistic for fixed eigenvalues $\Lambda$ and $D_0$. Another way to see the dependence on the eigenvectors is to rewrite $T$ as $T = 2n[\mathrm{tr}(\Lambda D_0) - \mathrm{tr}(\Lambda \tilde{V} D_0 \tilde{V}')]/\sigma^2$, where $\tilde{V} = V'U_0$. The second term in the brackets measures the angles between the eigenvectors $V$ and $U_0$ weighted by the eigenvalues. Here we can see how the multiplicities in $D_0$ play a crucial role. As the eigenvalues of $D_0$ get closer, the angles between the corresponding eigenvectors in $V$ and $U_0$ become harder to detect, and are unidentifiable if there are exact multiplicities.

It is also not surprising that $T$ in (31) does not depend on $\tau$. All the variability in $T$ is due to the variability of the MLE on $\mathcal{M}_{D_0}$, which is orthogonal to $\mathcal{I}_p$.

COROLLARY 4.2 [Case (S2)].    *Consider the test $H_0 : M \in \mathcal{M}_{D_0}$ vs. $H_A : M \notin \mathcal{M}_{D_0}$. Let $\bar{Y} = V\Lambda V'$ be an eigendecomposition of $\bar{Y}$ so that $\Lambda$ and $D_0$ are in decreasing order. Given $\sigma$ and $\tau$, the LLR test statistic is*

$$(32) \qquad T = n\|\Lambda - D_0\|^2_{\sigma^2,\tau} \xrightarrow[H_0]{n \to \infty} \chi^2\left(\frac{1}{2}\sum_{i=1}^k m_i(m_i+1)\right).$$

PROOF.    In the notation of Lemma 3.4, $\hat{M}_0 = VD_0V'$ by Theorem 4.1, and $\hat{M}_A = \bar{Y}$. Replacing these in (18a) gives (32). Under $H_0$, $\mathcal{M}_{D_0}$ has dimension $\dim(\mathcal{M}_{D_0}) = q - \sum_{i=1}^k m_i(m_i+1)/2$. Under $H_A$ the dimension is $q$, so the number of degrees of freedom for the test is $\sum_{i=1}^k m_i(m_i+1)/2$. $\square$

Since this is a test of whether $M_0$ has eigenvalues $D_0$ while treating the eigenvectors as nuisance parameters, it is not surprising that the test statistic (32) is equal to the normalized distance between $D_0$ and the eigenvalues of $\bar{Y}$.

### 4.4. *Curved submanifolds [Case (S3)].*

We have the embedded submanifold $\mathcal{M}_{m_1,\ldots,m_k}$ given by (5), the set of matrices with unspecified $k \le p$ distinct eigenvalues $\tilde{d}_1, \ldots, \tilde{d}_k$ and corresponding multiplicities $m_1, \ldots, m_k$. The eigenvalues are assumed to have a known order. For example, a prolate matrix ($d_1 > d_2 = d_3$) is different from an oblate matrix ($d_1 = d_2 > d_3$), even though both have the same number of multiplicities. The dimension of the set $\mathcal{M}_{m_1,\dots,m_k}$ is the dimension of the set defined in Lemma 4.1 plus $k$ for the free eigenvalues, in total

$$(33) \qquad \dim(\mathcal{M}_{m_1,\dots,m_k}) = k + q - \tfrac{1}{2}\sum_{i=1}^{k} m_i(m_i + 1).$$

The MLE of $M$ is given by Theorem 4.2 below. For this we need the following definition.

DEFINITION 4.1. Define the cumulative multiplicities $e_0 = 0$, $e_j = \sum_{i=1}^{j} m_i$, $j = 1,\dots,k$, so that $e_k = \sum_{i=1}^{k} m_i = p$. We define the *block average* of a diagonal matrix $\Lambda$ according to the ordered sequence of multiplicities $m_1,\dots,m_k$, denoted $\mathrm{blk}_{m_1,\dots,m_k}(\Lambda) \in \mathcal{D}_p$, as the block diagonal matrix formed by partitioning $\Lambda$ into diagonal blocks of sizes $m_1,\dots,m_k$ in that order and replacing the diagonal entries in each block by the average of the diagonal entries in that block. Specifically, for $j = 1,\dots,k$, the $j$th diagonal block of $\mathrm{blk}_{m_1,\dots,m_k}(\Lambda)$ is $(1/m_j)\sum_{i=e_{j-1}+1}^{e_j} \lambda_i I_{m_j}$. Below we use the shorthand notation $\mathrm{blk}(\Lambda)$ when the multiplicities are understood from the context.

THEOREM 4.2. Let $\mathcal{M}_{m_1,\dots,m_k}$ be given by (5) and let $\bar{Y} = V\Lambda V'$ be an eigendecomposition with eigenvalues in decreasing order. The MLE of $M$ is $\hat{M} = \hat{U}\hat{D}\hat{U}'$ where $\hat{U} = VQ$ and $Q$ is block diagonal orthogonal with orthogonal diagonal blocks $Q_1,\dots,Q_k$ of sizes $m_1,\dots,m_k$ and $\hat{D} = \mathrm{blk}_{m_1,\dots,m_k}(\Lambda)$ given by Definition 4.1.

PROOF. The set $\mathcal{M}_{m_1,\dots,m_k}$ is parallel to $\mathcal{I}_p$. By Theorem 3.1, the MLE minimizes $g(U,D) = \mathrm{tr}[(\bar{Y} - UDU')^2]$. Given $D$, the minimization with respect to $U$ is the same as in Theorem 4.1. Notice that the result does not depend on the actual entries of $D$ but only on their multiplicities. To find the MLE of $D$, choose $Q = I_p$ so that $\hat{U} = V$, where the eigenvectors correspond to eigenvalues in decreasing order. Then

$$g(\hat{U}, D) = \mathrm{tr}[(\bar{Y} - VDV')^2] = \mathrm{tr}[(\Lambda - D)^2] = \sum_{j=1}^{k}\sum_{i=e_{j-1}+1}^{e_j} (\lambda_i - \tilde{d}_j)^2.$$

Taking the derivative with respect to any particular eigenvalue $\tilde{d}_j$ and equating to zero gives $\hat{d}_j = (1/m_j)\sum_{i=e_{j-1}+1}^{e_j} \lambda_i$, which gets repeated over the block of size $m_j$. □

The solution to the test of whether $M$ is isotropic or partially isotropic with particular multiplicities of its eigenvalues is given by the next corollary. The test statistic is the distance between the eigenvalues of $\bar{Y}$ and the estimated eigenvalues under the assumption of the particular multiplicities.

COROLLARY 4.3 [Case (S3)].    *Consider the test* $H_0 : M \in \mathcal{M}_{m_1,\ldots,m_k}$ *vs.* $H_A : M \notin \mathcal{M}_{m_1,\ldots,m_k}$. *Given* $\sigma$ *and* $\tau$, *the LLR test statistic is*

$$(34) \qquad T = \frac{n}{\sigma^2} \operatorname{tr}[(\Lambda - \operatorname{blk}(\Lambda))^2] \xrightarrow[H_0]{n \to \infty} \chi^2 \left( \frac{1}{2} \sum_{i=1}^{k} m_i(m_i + 1) - k \right).$$

PROOF.    Replacing the MLE from Theorem 4.2 for $\hat{M}_0$ in (18a) and $\hat{M}_A = \bar{Y}$ gives $T = n \|\Lambda - \operatorname{blk}(\Lambda)\|^2_{\sigma^2, \tau}$. But the argument of the norm has trace zero, yielding (34). The number of degrees of freedom is equal to $q$ minus the dimension of $\mathcal{M}_{m_1,\ldots,m_k}$, which is given by (33).    $\square$

## 5. Inference for eigenvalues and eigenvectors in the two-sample problem.
As in the one-sample case, we assume $\sigma$ and $\tau$ are known. If not, adjustments for unknown $\sigma$ and $\tau$ can be made based on Lemma 3.6. The unrestricted case (A0) was covered in Example 3.2. We proceed with cases (S1) and (S2).

In order to solve the LLR statistics for these cases, we first need to find the MLEs of $M_1$ and $M_2$ in $\mathcal{M}_{2,D}$. The solution is given in Theorem 5.1 below. In what follows, let $\bar{Y}_1 = V_1 \Lambda_1 V_1'$, $\bar{Y}_2 = V_2 \Lambda_2 V_2'$ and $\bar{Y} = V \Lambda V'$ be eigendecompositions, all with eigenvalues in decreasing order.

THEOREM 5.1.    *Let* $\mathcal{M}_{2,D}$ *be given by* (7) *and suppose the eigenvalues* $\tilde{d}_1 > \cdots > \tilde{d}_k$ *of* $D$ *have known multiplicities* $m_1, \ldots, m_k$. *The MLEs of* $M_1$ *and* $M_2$ *in* $\mathcal{M}_{2,D}$ *are* $\hat{M}_1 = \hat{U}_1 \hat{D} \hat{U}_1'$ *and* $\hat{M}_2 = \hat{U}_2 \hat{D} \hat{U}_2'$ *where:*

(i) $\hat{U}_1$ *and* $\hat{U}_2$ *are any matrices of the form* $\hat{U}_1 = V_1 Q_1$ *and* $\hat{U}_2 = V_2 Q_2$, *where* $Q_1$ *and* $Q_2$ *are orthogonal matrices such that* $Q_1 D Q_1' = D$ *and* $Q_2 D Q_2' = D$, *that is, block diagonal orthogonal with orthogonal blocks of sizes* $m_1, \ldots, m_k$.

(ii) *Let* $\bar{\Lambda} = (n_1 \Lambda_1 + n_2 \Lambda_2)/n$. *Then* $\hat{D} = \operatorname{blk}_{m_1,\ldots,m_k}(\bar{\Lambda})$ *given by Definition 4.1.*

Even though $D$ is unspecified, it turns out that the MLE, while it does not depend on the actual value of $D$, does depend on the multiplicities of the eigenvalues of $D$, because that affects the dimension of the set $\mathcal{M}_{2,D}$. Since this set is parallel to $(\mathcal{I}_p, \mathcal{I}_p)$, the MLE does not depend on $\sigma$ or $\tau$.

PROOF OF THEOREM 5.1.    (i) $\mathcal{M}_{2,D}$ is parallel to $(\mathcal{I}_p, \mathcal{I}_p)$. By Lemma 3.7 and Theorem 3.2, the MLE minimizes

$$(35) \qquad g(U_1, U_2, D) = n_1 \operatorname{tr}[(\bar{Y}_1 - U_1 D U_1')^2] + n_2 \operatorname{tr}[(\bar{Y}_2 - U_2 D U_2')^2].$$

Given $D$, each summand is minimized separately as in Theorem 4.1, yielding $\hat{U}_1$ and $\hat{U}_2$. Notice that the result does not depend on the actual entries of $D$ but only on their multiplicities.

(ii) The MLE of $D$ minimizes (35). For the MLEs of $U_1$ and $U_2$, choose $Q_1 = Q_2 = I_p$ so that $\hat{U}_1 = V_1$ and $\hat{U}_2 = V_2$, where these eigenvectors correspond to eigenvalues in decreasing order. Replacing these in (35) gives

$$g(\hat{U}_1, \hat{U}_2, D) = n_1 \operatorname{tr}[(\Lambda_1 - D)^2] + n_2 \operatorname{tr}[(\Lambda_2 - D)^2]$$

$$= n_1 \sum_{j=1}^{k} \sum_{i=e_{j-1}+1}^{e_j} (\lambda_{1,i} - \tilde{d}_j)^2 + n_2 \sum_{j=1}^{k} \sum_{i=e_{j-1}+1}^{e_j} (\lambda_{2,i} - \tilde{d}_j)^2.$$

Taking the derivative with respect to a particular eigenvalue $\tilde{d}_j$ and equating to zero gives

$$\hat{d}_j = \frac{1}{m_j} \sum_{i=e_{j-1}+1}^{e_j} \frac{n_1 \lambda_{1,i} + n_2 \lambda_{2,i}}{n}. \qquad \square$$

We can now solve the LLR statistics.

COROLLARY 5.1 [Case (S1)].  *Consider the test $H_0 : M \in \mathcal{M}_{2,D}$ vs. $H_A : M \notin \mathcal{M}_{2,D}$, where $\mathcal{M}_{2,D}$ is given by (7). Given $\sigma$ and $\tau$, the LLR test statistic is*

$$T = \frac{n_1 n_2}{n} \|\Lambda_1 - \Lambda_2\|_{\sigma^2, \tau}^2 + \frac{n}{\sigma^2} \operatorname{tr}[(\Lambda - \operatorname{blk}(\bar{\Lambda}))^2] \tag{36}$$

*and is asymptotically $\chi^2(\sum_{i=1}^{k} m_i(m_i + 1) - k)$ under $H_0$ as $n \to \infty$.*

PROOF.  Under $H_0$, the MLEs are given by Theorem 5.1. Under $H_A$, both the eigenvalues and eigenvectors are unrestricted, and so the MLEs are simply $\hat{M}_{1,A} = \bar{Y}_1$ and $\hat{M}_{2,A} = \bar{Y}_2$. By Lemma 3.7, the LLR test statistic is

$$T = n_1 \|V_1 \Lambda_1 V_1' - V_1 \operatorname{blk}(\bar{\Lambda}) V_1'\|_{\sigma^2, \tau}^2 + n_2 \|V_2 \Lambda_2 V_2' - V_2 \operatorname{blk}(\bar{\Lambda}) V_2'\|_{\sigma^2, \tau}^2$$

$$= n_1 \|\Lambda_1 - \operatorname{blk}(\bar{\Lambda})\|_{\sigma^2, \tau}^2 + n_2 \|\Lambda_2 - \operatorname{blk}(\bar{\Lambda})\|_{\sigma^2, \tau}^2$$

$$= \frac{n_1 n_2}{n} \|\Lambda_1 - \Lambda_2\|_{\sigma^2, \tau}^2 + n \|\bar{\Lambda} - \operatorname{blk}(\bar{\Lambda})\|_{\sigma^2, \tau}^2,$$

where we have applied the equivalent form (21b) of the norm. The expression simplifies to (36) by noting that the second summand has trace zero. Under $H_0$, by Lemma 4.1, the parameter set has dimension $\dim(\mathcal{M}_{2,D}) = k + 2(q - \sum_{i=1}^{k} m_i(m_i + 1)/2) = k + 2q - \sum_{i=1}^{k} m_i(m_i + 1)$. Under $H_A$, the parameter set has dimension $\dim(\mathcal{M}_A) = 2q$. By Lemma 3.7, $T$ is asymptotically $\chi^2$ with number of degrees of freedom $\dim(\mathcal{M}_A) - \dim(\mathcal{M}_0) = \sum_{i=1}^{k} m_i(m_i + 1) - k$ degrees of freedom.  $\square$

The first term in the test statistic (36) resembles the test statistic for the unrestricted case (Example 3.2), but it involves the eigenvalues only. The second term captures the effect of the multiplicities. If $D$ has no repeats (all the eigenvalues have multiplicity one), then the block average of $\bar{\Lambda}$ is the same as $\Lambda$ and the second term drops. The asymptotic distribution of the test statistic in this case is $\chi^2(p)$. On the other hand, if $D$ is fully isotropic (one eigenvalue of multiplicity $p$), we get a $\chi^2(2q-1)$.

COROLLARY 5.2 [Case (S2)].   *Consider the test* $H_0 : M_1 = M_2$ *vs.* $H_A : (M_1, M_2) \in \mathcal{M}_{2,D}$ *where* $\mathcal{M}_{2,D}$ *is given by* (7). *Under both* $H_0$ *and* $H_A$, $D$ *is assumed to have unknown eigenvalues* $\bar{d}_1 > \cdots > \bar{d}_k$ *with known multiplicities* $m_1, \ldots, m_k$. *Given* $\sigma$ *and* $\tau$, *the LLR test statistic is*

$$
\begin{aligned}
(37) \qquad T = {} & \frac{2n_1 n_2}{n\sigma^2}[\operatorname{tr}(\Lambda_1 \Lambda_2) - \operatorname{tr}(\bar{Y}_1 \bar{Y}_2)] \\
& + \frac{n}{\sigma^2}\operatorname{tr}[(\Lambda - \operatorname{blk}(\Lambda))^2 - (\bar{\Lambda} - \operatorname{blk}(\bar{\Lambda}))^2]
\end{aligned}
$$

*and is asymptotically* $\chi^2(q - \sum_{i=1}^{k} m_i(m_i+1)/2)$ *under* $H_0$ *as* $n \to \infty$.

PROOF.   Under $H_0$, the data can be regarded as one sample of size $n$. By Theorem 4.2, the MLEs are $\hat{M}_{1,0} = \hat{M}_{2,0} = V\operatorname{blk}(\Lambda)V'$. Under $H_A$, Theorem 5.1 gives $\hat{M}_{1,A} = V_1\operatorname{blk}(\bar{\Lambda})V_1'$ and $\hat{M}_{2,A} = V_2\operatorname{blk}(\bar{\Lambda})V_2'$. By Lemma 3.7, replacing these MLEs in (25) gives

$$
\begin{aligned}
T = {} & n\|\bar{Y} - V\operatorname{blk}(\Lambda)V'\|_{\sigma^2,\tau}^2 + \frac{n_1 n_2}{n}\|\Delta(\bar{Y})\|_{\sigma^2,\tau}^2 \\
& - n_1\|\Lambda_1 - \operatorname{blk}(\bar{\Lambda})\|_{\sigma^2,\tau}^2 - n_2\|\Lambda_2 - \operatorname{blk}(\bar{\Lambda})\|_{\sigma^2,\tau}^2 \\
= {} & n\|\Lambda - \operatorname{blk}(\Lambda)\|_{\sigma^2,\tau}^2 + \frac{n_1 n_2}{n}\|\bar{Y}_1 - \bar{Y}_2\|_{\sigma^2,\tau}^2 \\
& - n\|\bar{\Lambda} - \operatorname{blk}(\bar{\Lambda})\|_{\sigma^2,\tau}^2 - \frac{n_1 n_2}{n}\|\Lambda_1 - \Lambda_2\|_{\sigma^2,\tau}^2 \\
= {} & \frac{2n_1 n_2}{n}[\langle \Lambda_1, \Lambda_2 \rangle_{\sigma^2,\tau} - \langle \bar{Y}_1, \bar{Y}_2 \rangle_{\sigma^2,\tau}] \\
& + n\|\Lambda - \operatorname{blk}(\Lambda)\|_{\sigma^2,\tau}^2 - n\|\bar{\Lambda} - \operatorname{blk}(\bar{\Lambda})\|_{\sigma^2,\tau}^2,
\end{aligned}
$$

where we have used the form (21b) twice. Expression (37) follows by noting that $\operatorname{tr}(\Lambda_1)\operatorname{tr}(\Lambda_2) = \operatorname{tr}(\bar{Y}_1)\operatorname{tr}(\bar{Y}_2)$ and that the terms inside the norms have trace zero.

Under $H_0$, by Lemma 4.1, the parameter set has dimension $\dim(\mathcal{M}_0) = k + q - \sum_{i=1}^{k} m_i(m_1+1)/2$. Under $H_A$, the parameter set has dimension $\dim(\mathcal{M}) = k + 2(q - \sum_{i=1}^{k} m_i(m_i+1)/2) = k + 2q - \sum_{i=1}^{k} m_i(m_i+1)$. By Lemma 3.7, $T$ is asymptotically $\chi^2$ with number of degrees of freedom $\dim(\mathcal{M}_A) - \dim(\mathcal{M}_0) = q - \sum_{i=1}^{k} m_i(m_i+1)/2$ degrees of freedom.   □

The form of the test statistic (37) is interesting. Recall that (S2) is a test of whether $M_1$ and $M_2$ have the same eigenvectors when they are assumed to have the same eigenvalues. Consider the first bracket $\tilde{T}_1 = \operatorname{tr}(\Lambda_1 \Lambda_2) - \operatorname{tr}(\bar{Y}_1 \bar{Y}_2)$. The behavior of this term is similar to (31). If the eigenvectors in $V_1$ and $V_2$ are the same, then $\tilde{T}_1$ is zero. As the angles between $V_1$ and $V_2$ increase, then the inner product between $\bar{Y}_1$ and $\bar{Y}_2$ decreases, which increases the value of $T_1$ for fixed eigenvalues $\Lambda_1$ and $\Lambda_2$. Another way to see the dependence on the eigenvectors is to rewrite $\tilde{T}_1$ as $\tilde{T}_1 = \operatorname{tr}(\Lambda_1 \Lambda_2) - \operatorname{tr}(\Lambda_1 \tilde{V} \Lambda_2 \tilde{V}')$, where $\tilde{V} = V_1' V_2$. The second term measures the angles between the eigenvectors $V_1$ and $V_2$ weighted by the eigenvalues. Here we see how the multiplicities of the common unknown eigenvalues $D$ play a crucial role. As the eigenvalues of $D$ get closer, the angles between the corresponding eigenvectors in $V_1$ and $V_2$ become harder to detect, and are unidentifiable if there are exact multiplicities. The multiplicities also play a role in the last two terms in (37). If the eigenvectors $V_1$ and $V_2$ are equal (up to sign), then the difference between those two terms is zero. As the angles between $V_1$ and $V_2$ increase, the eigenvalues $\Lambda$ of the pooled average decrease with respect to the average of the eigenvalues $\bar{\Lambda}$, which also increases the value of $T$.

As mentioned earlier in Section 2.2, cases (S3) and (S4) are difficult because they involve the set $\mathcal{M}_{2,U}$ given by (8), whose MLE has no closed-form solution. The MLE for the common matrix $U \in \mathcal{O}_p$ in this case could potentially be found using numerical optimization techniques such as the ones described by Edelman, Arias and Smith [10]. A useful simplification here is that the set $\mathcal{M}_{2,U}$ is parallel to $\mathcal{I}_p \times \mathcal{I}_p$ in the sense of Definition 3.2, so by Theorem 3.2, the MLE given $\sigma$ and $\tau$ can be solved assuming $\sigma^2 = 1$ and $\tau = 0$.

**6. Summary and discussion.** In this article we have derived MLEs and LLR tests for the mean parameter $M$ of the orthogonally invariant symmetric-matrix-variate normal distribution. This has been carried out for many families of interest regarding eigenvalues and eigenvectors of $M$, both in the one-sample and two-sample settings. The parameter sets involved have been affine subspaces and orthogonally invariant embedded submanifolds of $\mathcal{S}_p$ and $\mathcal{S}_p \times \mathcal{S}_p$, in addition to one case in which the parameter set is a polyhedral convex cone in $\mathcal{S}_p$. In these parameter sets, the geometry of the sets and the multiplicity pattern of the eigenvalues have played crucial roles. Despite some of the classical characteristics of this problem, conceivably it was not attempted before because there was no data that required it. The emergence of new data such as DTI and polarized CBR has made this need come to life.

As pointed out by Theorems 3.1 and 3.2, the main advantage of adopting an orthogonally invariant covariance structure has been analytical, allowing the derivation of explicit formulas for the MLEs and LLRs. We have identified a class of sets in $\mathcal{S}_p$ and another in $\mathcal{S}_p \times \mathcal{S}_p$ where the MLEs are easy to obtain because they do not depend on the covariance parameters. Furthermore, Theorems 3.1 and 3.2 show that the orthogonally invariant covariance is the most general covariance structure with that property within those classes of sets. The explicit formulas presented in Propositions 3.1, 4.1 and 4.2, in Theorems 4.1, 4.2 and 5.1 and in Corollaries 4.1, 4.2, 4.3, 5.1 and 5.2 have helped gain analytical insights. But they are also useful from a computational point of view. Both in DTI and CBR data, images contain hundreds of thousands of pixels, with one random symmetric matrix at each pixel. For computations to be efficient, it is extremely useful to have explicit formulas that can be easily evaluated at each pixel.

For real data, the normality assumption can be checked for any particular data set using any test of multivariate normality applied to the data vectors $\text{vecd}(Y_i)$, $i = 1, \ldots, n$, defined in Section 3. While the orthogonally invariant covariance structure has been shown to be appropriate for DTI (Basser and Pajevic [3]), this assumption can also be checked using the LLR test for orthogonal invariance given in Proposition 3.1. More ambitious data analysts may not want to make such assumptions on the covariance. In those cases, least squares solutions may be found using numerical optimization methods. Another option is to use the statistics derived in this paper but to adjust their distributions when the true covariance is not orthogonally invariant. This would be akin to using the standard least squares solution in a linear regression problem with correlated errors and then adjusting the standard errors of the estimates. Moreover, the MLEs in this paper can still serve as least squares estimators when the distribution of the data is not Gaussian. These lines of research are left for future work.

**Acknowledgments.** This article is based on a chapter of the first author's Ph.D. thesis [30], written at Stanford University under the supervision of Professors Bradley Efron and Jonathan Taylor. The authors wish to thank Bradley Efron for his valuable guidance, Donald Richards for useful discussions and the referees for their careful revisions and helpful remarks.

## REFERENCES

[1] ARSIGNY, V., FILLARD, P., PENNEC, X. and AYACHE, N. (2005). Fast and simple calculus on tensors in the log-Euclidean framework. *MICCAI 2005. Lecture Notes in Comput. Sci.* **3749** 115–122.

[2] ARSIGNY, V., FILLARD, P., PENNEC, X. and AYACHE, N. (2007). Geometric means in a novel vector space structure on symmetric positive definite matrices. *SIAM. J. Matrix Anal. Appl.* **29** 328–347.

[3] BASSER, P. J. and PAJEVIC, S. (2003). A normal distribution for tensor-valued random variables: applications to diffusion tensor MRI. *IEEE Trans. Med. Imaging* **22** 785–794.

[4] Basser, P. J. and Pierpaoli, C. (1996). Microstructural and physiological features of tissues elucidated by quantitative–diffusion–tensor MRI. *J. Magn. Reson. B* **111** 209–219.

[5] Chang, T. (1986). Spherical regression. *Ann. Statist.* **14** 907–924. MR0856797

[6] Chernoff, H. (1954). On the distribution of the likelihood ratio. *Ann. Math. Statist.* **25** 573–578. MR0065087

[7] Chikuse, Y. (2003). *Statistics on Special Manifolds*. Springer, New York. MR1960435

[8] Drton, M. (2008). Likelihood ratio tests and singularities. *Ann. Statist.* To appear.

[9] Dykstra, R. L. (1983). An algorithm for restricted least squares regression. *J. Amer. Statist. Assoc.* **78** 837–842. MR0727568

[10] Edelman, A., Arias, T. A. and Smith, S. T. (1998). The geometry of algorithms with orthogonality constraints. *SIAM J. Matrix Anal. Appl.* **20** 303–353. MR1646856

[11] Efron, B. (1978). The geometry of exponential families. *Ann. Statist.* **6** 362–376. MR0471152

[12] Fang, K.-T. and Zhang, Y.-T. (1990). *Generalized Multivariate Analysis*. Springer, Berlin. MR1079542

[13] Fletcher, P. T. and Joshi, S. (2007). Riemannian geometry for the statistical analysis of diffusion tensor data. *Signal Processing* **87** 250–262.

[14] Gupta, A. K. and Nagar, D. K. (2000). *Matrix Variate Distributions*. Chapman and Hall/CRC Publisher, Boca Raton, FL. MR1738933

[15] Hu, W. and White, M. (1997). A CMB polarization primer. *New Astronomy* **2** 323–344.

[16] James, A. T. (1976). Special functions of matrix and single argument in statistics. In *Theory and Applications of Special Functions* (R. A. Askey, ed.) 497–520. Academic Press, New York. MR0402145

[17] Kogut, A., Spergel, D. N., Barnes, C., Bennett, C. L., Halpern, M., Hinshaw, G., Jarosik, N., Limon, M., Meyer, S. S., Page, L., Tucker, G. S., Wollack, E. and Wright, E. L. (2003). Wilkinson microwave anisotropy probe (WMAP) first year observations: TE polarization. *Astrophysical J. Supplement Series* **148** 161–173.

[18] Lang, S. (1999). *Fundamentals of Differential Geometry*. Springer, New York. MR1666820

[19] Lawson, C. L. and Hanson, B. J. (1974). *Solving Least Squares Problems*. Prentice-Hall Inc., Englewood Cliffs, NJ. MR0366019

[20] LeBihan, D., Mangin, J.-F., Poupon, C., Clark, C. A., Pappata, S., Molko, N. and Chabriat, H. (2001). Diffusion tensor imaging: Concepts and applications. *J. Magn. Reson. Imaging* **13** 534–546.

[21] Lehman, E. L. (1997). *Testing Statistical Hypotheses*, 2nd ed. Springer, New York. MR1481711

[22] Mallows, C. L. (1961). Latent vectors of random symmetric matrices. *Biometrika* **48** 133–149. MR0131312

[23] Mardia, K. V., Kent, J. T. and Bibby, J. M. (1979). *Multivariate Analysis*. Academic Press, San Diego, CA. MR0560319

[24] Mehta, M. L. (1991). *Random Matrices*, 2nd ed. Academic Press, San Diego, CA. MR1083764

[25] Moakher, M. (2002). Means and averaging in the group of rotations. *SIAM J. Matrix Anal. Appl.* **24** 1–16. MR1920548

[26] Pajevic, S. and Basser, P. J. (2003). Parametric and nonparametric statistical analysis of DT-MRI data. *J. Magn. Reson.* **161** 1–14.

[27] RAUBERTAS, R. R. (2006). Pool-adjacent-violators algorithm. In *Encyclopedia of Statistical Sciences*. Wiley, New York.

[28] ROBERTSON, T. and WEGMAN, E. J. (1978). Likelihood ratio tests for order restructions in exponential families. *Anns. Statist.* **6** 485–505.

[29] SCHEFFÉ, H. (1970). Practical solutions to the Behrens–Fisher problem. *J. Amer. Statist. Assoc.* **65** 1501–1508. MR0273732

[30] SCHWARTZMAN, A. (2006). *Random ellipsoids and false discovery rates: statistics for diffusion tensor imaging data*. Ph.D. dissertation, Stanford Univ.

[31] SELF, S. G. and LIANG, K.-Y. (1987). Asymptotic properties of maximum likelihood estimators and likelihood ratio tests under nonstandard conditions. *J. Amer. Statist. Assoc.* **82** 605–610. MR0898365

[32] WALD, A. (1949). Note on the consistency of the maximum likelihood estimate. *Ann. Math. Statist.* **20** 595–601. MR0032169

[33] WHITCHER, B., WISCO, J. J., HADJIKHANI, N. and TUCH, D. S. (2007). Statistical group comparison of diffusion tensors via multivariate hypothesis testing. *Magn. Reson. Med.* **57** 1065–1074.

[34] ZHU, H., ZHANG, H., IBRAHIM, J. G. and PETERSON, B. S. (2007). Statistical analysis of diffusion in diffusion-weighted magnetic resonance imaging data. *J. Amer. Statist. Assoc.* **102** 1085–1102.

A. SCHWARTZMAN
DEPARTMENT OF BIOSTATISTICS
HARVARD SCHOOL OF PUBLIC HEALTH
BOSTON, MASSACHUSETTS 02115
AND
DEPARTMENT OF BIOSTATISTICS
    AND COMPUTATIONAL BIOLOGY
DANA-FARBER CANCER INSTITUTE
44 BINNEY STREET, CLS-11007
BOSTON, MASSACHUSETTS 02115
USA
E-MAIL: armins@hsph.harvard.edu

W. F. MASCARENHAS
DEPARTAMENTO DE CIÊNCIA DA COMPUTAÇÃO
INSTITUTO DE MATEMÁTICA E ESTATÍSTICA
UNIVERSIDADE DE SÃO PAULO
RUA DO MATÃO 1010
05508-090 SÃO PAULO
BRAZIL
E-MAIL: walterfm@ime.usp.br

J. E. TAYLOR
DEPARTMENT OF STATISTICS
STANFORD UNIVERSITY
390 SERRA MALL
STANFORD, CALIFORNIA 94305
USA
E-MAIL: jonathan.taylor@stanford.edu